\documentclass[11pt,dvips,leqno]{article}

\usepackage{amsmath}
\usepackage{xspace}
\usepackage{amssymb,latexsym}
\usepackage{graphicx}
\usepackage{psfrag}
\usepackage{colordvi}
\usepackage{color}

\def\nothanks#1{}

\date{4 March 2003}

\newcommand{\distr}{\sim}

\def\Box{\square}

\newtheorem{theorem}{Theorem}

\newtheorem{claim}[theorem]{Claim}
\newtheorem{conjecture}[theorem]{Conjecture}
\newtheorem{remark}[theorem]{Remark}

\renewcommand{\leq}{\leqslant}
\renewcommand{\geq}{\geqslant}

\def\ra{\rightarrow}

\newcommand{\E}{\mathbb{E}}
\newcommand{\pr}{\mathbb{P}}
\renewcommand{\Re}{\mathbb{R}}

\newcommand{\given}{\;|\;}
\newcommand{\ve}{\varepsilon}

\newcommand{\Opt}{\ensuremath{\operatorname{Opt}}\xspace}

\newenvironment{proof}%
  {\smallbreak\noindent \textbf{Proof.\ }}%
  {\hspace*{1em} \hfill $\Box$\vspace{0.3cm}}

\newcommand{\eps}{\varepsilon}

\newcommand{\whp}{w.h.p.\xspace}
\newcommand{\as}{a.s.\xspace}

\renewcommand{\a}{\alpha}
\newcommand{\boll}{Bollob\'{a}s\xspace}
\newcommand{\bw}{\bar{w}}
\newcommand{\bu}{\bar{u}}
\newcommand{\bv}{\bar{v}}
\newcommand{\cn}{\lfloor cn \rfloor}

\newcommand{\csp}{\textnormal{\textsc{csp}}\xspace}
\newcommand{\cut}{\operatorname{cut}}
\newcommand{\D}{\Delta}
\renewcommand{\d}{\delta}
\newcommand{\e}{\varepsilon}
\newcommand{\eqdef}{\doteq}
\newcommand{\fc}{f_{\cut}}
\newcommand{\nogiant}{\operatorname{nogiant}}
\newcommand{\fnogiant}{f_{\nogiant}}
\newcommand{\foa}{f_{\textsc{O-I}}}
\newcommand{\fob}{f_{\textsc{O-II}}}
\renewcommand{\l}{\lambda}
\newcommand{\logis}[1]{{#1} \ln{#1}}
\newcommand{\maxsc}{\textnormal{\textsc{max}}\xspace}
\newcommand{\minsc}{\textnormal{\textsc{min}}\xspace}
\newcommand{\maxcsp}{\maxsc \csp}
\newcommand{\maxksat}{\textsc{max} $k$-\textsc{sat}\xspace}
\newcommand{\maxsat}{\textnormal{\textsc{max sat}}\xspace}
\newcommand{\maxtsat}{\textnormal{\textsc{max 2-sat}}\xspace}
\newcommand{\maxcut}{\textnormal{\textsc{max cut}}\xspace}
\newcommand{\mgfss}{\maxsc giant-free spanning subgraph\xspace}

\newcommand{\MGF}{\maxsc giant-free\xspace}
\newcommand{\np}{\textsc{NP}}

\newcommand{\onlinea}{\textsc{Online-Lazy}\xspace}
\newcommand{\problema}{\textsc{Online~I}\xspace}
\newcommand{\problemb}{\textsc{Online~II}\xspace}
\newcommand{\rhostar}{\rho^\star}
\newcommand{\set}[1]{\left\{#1\right\}}
\newcommand{\F}{\mathcal{F}}
\newcommand{\const}{\mathop{\mathrm{const}}}
\newcommand{\Abar}{\bar{A}}

\newcommand{\true}{\text{True}\xspace}
\newcommand{\false}{\text{False}\xspace}
\newcommand{\sat}{\textnormal{\textsc{sat}}\xspace}
\newcommand{\Var}{\operatorname{Var}}
\newcommand{\w}{w}
\newcommand{\Xbar}{\bar{X}}
\newcommand{\Xvec}{\vec{X}}

\newcommand{\geqa}{\gtrsim}
\newcommand{\leqa}{\lesssim}
\renewcommand{\asymp}{\simeq}
\newcommand{\xor}{\textsc{xor}\xspace}
\newcommand{\fmath}[1]{${#1}$}

\newcommand{\sleq}{\preceq}
\newcommand{\one}{\mathbf{1}}
\renewcommand{\d}{\delta}
\renewcommand{\b}{\beta}
\renewcommand{\a}{\alpha}
\renewcommand{\t}{\tau}

\newcommand{\gnm}{\ensuremath{G(n,m)}\xspace}
\newcommand{\gnp}{\ensuremath{G(n,p)}\xspace}

\newcommand{\expfn}[1]{\exp\left({#1}\right)}
\newcommand{\Be}{\operatorname{Be}}
\newcommand{\Po}{\operatorname{Po}}
\newcommand{\clausein}{\subseteq}
\newcommand{\subsubsubsection}[1]{\subsubsection*{#1}}

\let\safelabel=\label
\newcommand{\dpf}{{\framebox{\small pf $\ra$ appendix}\hspace*{.2in}}}
\newlength{\canwidth}
\long\def\pad#1#2{\ifcat{#2}{ }{ #1 }\else{ #1#2}\fi}
\newcommand{\boldref}[1]{{{\global\let\label=\ignore}\textbf{\ref{#1}.} }}
\long\def\can#1#2{\global\long\def#1{#2}}
\newcommand{\copycan}[1]{{{%
    \vspace{1em}
    \let\label=\boldref
    \renewcommand{\dpf}{{}}
    \let\footnote=\ignore
    \renewenvironment{theorem}{\noindent
      \textbf{Theorem} \begin{em}}{\end{em}}
    \par \noindent {#1} \par%
    \global\let\label=\safelabel
}}}
\newcommand{\uncan}[1]{{#1}}
\newcommand{\pf}[1]{\vspace{1.0em} \noindent \textit{#1}}
\newcommand{\thmpf}[2]{\vspace{1.0em}
 \noindent \textsc{Theorem~\ref{#1}: \textit{#2}}}

\def\qed{\hspace*{1em} \hfill $\Box$}

\begin{document}
\makeatletter
\title{Random \maxsat, Random \maxcut, \\ and Their Phase Transitions}
\author{%
Don Coppersmith\thanks{Department of Mathematical Sciences,
IBM T.J.~Watson Research Center, Yorktown Heights NY 10598, USA.
\hbox{e-mail}~{\small\texttt{\{copper,gamarnik,sorkin\}@watson.ibm.com}}}
\and
David Gamarnik$^\ast$\nothanks{Department of Mathematical Sciences,
IBM T.J.~Watson Research Center, Yorktown Heights NY 10598, USA.
\hbox{e-mail}~{\small\texttt{gamarnik@watson.ibm.com}}}
\and
Mohammad Hajiaghayi\thanks{Department of Mathematics,
M.I.T., Cambridge MA 02139, USA.
\hbox{e-mail}~{\small\texttt{hajiagha@math.mit.edu}}}
\and
Gregory B. Sorkin$^\ast$\nothanks{Department of Mathematical Sciences,
IBM T.J.~Watson Research Center, Yorktown Heights NY 10598, USA.
\hbox{e-mail}~{\small\texttt{sorkin@watson.ibm.com}}}
}

\maketitle
\makeatother

\begin{abstract}
Given a 2-\sat formula $F$ consisting of $n$ variables and $\cn$
random clauses, what is the largest number of clauses $\max F$
satisfiable by a single assignment of the variables? We bound the
answer away from the trivial bounds of $\frac34 cn$ and $cn$. We
prove that for $c<1$, the expected number of clauses satisfiable
is $\cn-\Theta(1/n)$; for large $c$, it is $(\frac34 c +
\Theta(\sqrt{c}))n$; for $c = 1+\e$, it is at least
$(1+\e-O(\e^3))n$ and at most $(1+\e-\Omega(\e^3/\ln \e))n$; and
in the ``scaling window'' $c= 1+\Theta(n^{-1/3})$, it is
$cn-\Theta(1)$. In particular, just as the decision problem undergoes a
phase transition, our optimization problem also undergoes a phase
transition at the same critical value $c=1$.

Nearly all of our results are established without reference to
the analogous propositions for decision 2-\sat, and
as a byproduct we reproduce many of those results,
including much of what is known about the 2-\sat scaling window.

We consider ``online'' versions of \maxtsat, and show that for
one version, the obvious greedy algorithm is optimal.

We can extend only our simplest \maxtsat results to \maxsc $k$-\sat,
but we conjecture a ``\maxsc $k$-\sat limiting function conjecture''
analogous to the folklore satisfiability threshold conjecture,
but open even for $k=2$.
Neither conjecture immediately implies the other, but it is natural
to further conjecture a connection between them.

Finally, for random \maxcut
(the size of a maximum cut in a sparse random graph)
we prove analogous results.
\end{abstract}

\maketitle
\tableofcontents
\pagestyle{plain}

\section{Introduction}

In this paper, we consider random instances of \maxtsat,
\maxksat, and \maxcut.
Just as random instances of the decision problem 2-\sat
show a phase transition from almost-sure satisfiability to
almost-sure unsatisfiability as the instance ``density''
increases above~1,
so the maximization problem shows a transition at the
same point, with the expected number of clauses \emph{not} satisfied
by an optimal solution quickly
changing
from $\Theta(1/n)$ to $\Theta(n)$.
\textsc{Max cut} experiences a similar phase transition:
as a random graph's edge density crosses above $1/n$,
the number of edges \emph{not} cut in an optimal cut
suddenly changes from $\Theta(1)$ to $\Theta(n)$.

Our methods are well established ones:
the first-moment method for upper bounds;
algorithmic analysis including the differential-equation method
for lower bounds;
and some more sophisticated arguments for the analysis
of the scaling window.
The interest of the work lies in the simplicity of the methods,
and in the results.
The questions we ask seem very natural, and the
answers obtained for \maxtsat and \maxcut are happily neat,
and fairly comprehensive.

A preliminary version of this paper appeared as~\cite{SODA03max}.

\subsection{Outlook}
Beyond our particular results for \maxtsat and \maxcut,
we hope to spark further work on
phase transitions in random instances of other optimization problems,
in particular of \maxsc \csp{s} (constraint satisfaction problems).
Random instances of optimization problems have been studied
extensively --- some that come to mind are
the travelling salesman problem, minimum spanning tree,
minimum assignment, minimum bisection,
minimum coloring, and maximum clique ---
but little
has been said about \emph{phase transitions}
in such cases, and indeed many of the examples do not even have
a natural parameter whose continuous variation could give rise
to a phase transition.

Many problems, including all \csp{s}, have natural
decision and optimization versions:
one can ask
whether a graph is $k$-colorable,
or ask for the minimum number of colors it requires.
We suggest that in a random setting,
the optimization version is quite as interesting as the decision version.
Furthermore, optimization problems may plausibly be easier to analyze
than decision problems because the quantities of interest vary
more smoothly.
In fact, a recent triumph in the analysis of a decision problem,
the characterization of the ``scaling window'' for 2-\sat,
used as a smoothed quantity
the size of the ``spine'' of a formula~\cite{BBCKW01}.
A way to view our \maxtsat results is that instead of
taking the size of the spine as our ``order parameter'',
we take the size of a maximum satisfiable subformula.
This seems comparably tractable
(we reproduce the result of \cite{BBCKW01} incompletely,
but more easily),
and arguably more natural.
Generally, when a decision problem has an optimization analog,
the value of the optimum is both interesting in its own right,
and, we suggest, an obvious candidate order parameter for studying the
decision problem.

\subsection{Problem and motivations}
Let $F$ be a $k$-\sat formula with $n$ variables $X_1,\ldots,X_n$.
An ``assignment'' of these variables consists of setting each $X_i$
to either 1 (True) or 0 (False); we may write an assignment as
a vector $\Xvec \in \{0,1\}^n$.
$k$-\sat is well understood.  In particular, it is a canonical
\np-hard problem to determine if a given formula $F$ is
satisfiable or not, except for $k=2$ when this decision problem
is solvable in essentially linear time.

Random instances of $k$-\sat have recently received wide attention.
Let $\F(n,m)$ denote the set of all formulas with $n$ variables
and $m$ clauses,
where each clause is proper
(consisting of $k$ distinct variables, each
of which may be complemented or not),
and clauses may be repeated.
Let $F \in \F$ be chosen uniformly at random;
this is equivalent to choosing $m$ clauses uniformly at random,
with replacement,
from the $2^k \binom{n}{k}$ possible clauses.

The model is generally parametrized as $F \in \F(n,cn)$
for various ``densities'' $c$, and
the state of knowledge is summarized thus.
The 2-\sat case is well understood:
for $c<1$, $F$ is almost surely satisfiable
(a.s.\ in the limit $n \ra \infty$),
and for $c>1$, $F$ is \as unsatisfiable%
\cite{CR92,GoTU,vega2sat}.
Recently, the ``scaling window'' $c=1 \pm \Theta(n^{-1/3})$
has also been analyzed~\cite{BBCKW01}.
For $k$-\sat, much less is known.
For 3-\sat, for instance, it is known that
for $c<3.42$, $F$ is \as satisfiable~\cite{KKL02}
and for $c>4.6$, $F$ is \as unsatisfiable~\cite{JSV}.
It is only conjectured, though, that for $k=3$
(and for all $k$) the situation is similar to that for $k=2$.
\begin{conjecture}[Satisfiability Threshold Conjecture]
\label{conj:thresh}
For each $k$ there exists a threshold density $c_k$,
such that for any positive $\e$,
for all $c<c_k-\e$, a random formula $F$ is \as satisfiable,
and for all $c>c_k+\e$, $F$ is \as unsatisfiable.
\end{conjecture}
For large values of $k$, although the question
of a threshold remains open,
satisfiability and unsatisfiability density bounds
are asymptotically equal, as shown by an analysis
in~\cite{AchMoore} and refined in~\cite{AchPeres}.
The closest result to the satisfiability conjecture is
a theorem of Friedgut~\cite{Frie} proving similar thresholds,
but leaving open the possibility that (for a given~$k$),
each $n$ may have its own threshold,
and that these may not converge to a limit.
\begin{theorem}[Friedgut]
\label{thm:friedgut}
For each $k$ there exists a threshold density function $c_k(n)$,
such that for any positive $\e$,
as $n \ra \infty$,
for all $c<c_k-\e$, a random formula $F$ is \as satisfiable,
and for all $c>c_k+\e$, $F$ is \as unsatisfiable.
\end{theorem}

Having briefly surveyed \emph{random} $k$-\sat,
let us similarly consider \emph{max} $k$-sat.
For a given formula $F$,
let $F(\Xvec)$ be the number of clauses satisfied by~$\Xvec$.
The problem \maxtsat asks for $\max F \eqdef \max_{\Xvec} F(\Xvec)$,
i.e., the maximum,
over all assignments~$\Xvec$,
of the size (number of clauses)
of a maximum satisfiable subformula of~$F$.

In the maximization setting, even 2-\sat is interesting.
\maxtsat is \np-hard to solve exactly,
and it is even NP-hard to approximate $\max F$
to within a factor of $21/22$~\cite{Hastad97}.
On the other hand, a $3/4$-approximation is trivial:
a random assignment satisfies an expected $3/4$ths of the
clauses, and a derandomized algorithm is simple
(our algorithm used to prove the lower bound for
Theorem~\ref{MainHighDens} can serve).
The best known approximation ratio achievable in
polynomial time is~0.940~\cite{LLZ02}.
For arbitrary 3-\sat formulas $F$,
in polynomial time,
$\max F$ can be approximated to within a factor of
$7/8$~\cite{KZ97},
but no better (unless P$=$NP)~\cite{Hastad97}.

Although both randomized and maximization versions of $k$-\sat
are thus well studied, we are aware of no work on
random \maxsat, nor other random \maxsc or \minsc
constraint satisfaction problems (\csp{s}).
These problems seem very natural, and
answers to even the simplest questions are not
obvious at first blush:
For a random 2-\sat formula $F(n, c n)$ with $c>1$,
which is \as unsatisfiable, can we perhaps \whp satisfy
all but a single clause?

These questions have elegant answers; we will show for example
that random \maxtsat has a phase structure analogous
to the decision problem's.
And there is a hope that the maximization problems may help
in understanding the decision problems.
For 2-\sat this hope is borne out to a degree by our
Theorem~\ref{MainScalingWindow}.
While our results for $k>2$ are very limited
(see Theorem~\ref{ksat1}),
Conjectures~\ref{conjecture:limit} and~\ref{conj2} link the open questions
for the maximization and decision thresholds for random satisfiability.
At this point we cannot guess the comparative difficulties of resolving
the satisfiability threshold conjecture, its maximization analog,
or the conjectured link between them.

Our study of random \maxtsat and random \maxcut was also motivated
by recent work on ``avoiding a giant component'';
we will discuss this in section~\ref{sec:cut}.

We consider several aspects of random \maxtsat and random \maxcut.
We also extend the easiest results to
arbitrary {\csp}s (constraint satisfaction problems).

We will give a second motivation for considering problems of this
sort when we take up \maxcut, in section~\ref{sec:cut}.

\section{Notation and model}

We write $F(n,m)$ to denote a random 2-\sat formula
on $n$ variables, with $m$ clauses.
Typically we will fix a constant
$c$ and consider $F(n,\cn)$; where it does not matter we will
often write $cn$ in lieu of $\cn$ and we often omit the notation
$\lfloor\cdot\rfloor$ in other instances.
For any formula $F$, define $\max F$ to be the size of
a largest satisfiable subformula of~$F$.
Our focus is the functional behavior of $\max F$.

Similarly, we write \gnm for a random graph on $n$ vertices
with $m$ edges.
For any graph $G$, let $\Xvec$ describe a partition of the vertices,
and let $\cut(G,\Xvec)$ be the number of edges having one vertex in
each part of the partition.
Define $\max \cut(G) \eqdef \max_{\Xvec} \cut(G,\Xvec)$,
and $\fc(n,m) \eqdef \E( \max \cut(G(n,m)) )$.

We use standard asymptotic and ``order'' notation,
so for example $f(n) \asymp g(n)$ means $f(n)/g(n) \ra 1$
as $n \ra \infty$,
and $f(n) = o(n)$ means $f(n)/n \ra 0$.
We will also write $f(n) \leqa g(n)$ to indicate that
$f$ is less than or equal to $g$ \emph{asymptotically}
--- $\limsup f(n)/g(n) \leq 1$ ---
though it may be that $f(n) > g(n)$ even for arbitrarily large values of~$n$.
Asymptotic results involving two variables, for example concerning
2-\sat formulae on $n$ variables with $cn$ clauses, with $c$ large
(or $(1+\e)n$ clauses with $\e$ small) should always be
interpreted as taking the limit in $n$ second; thus ``for any
desired error bound there exists a $c_0$, such that for all
$c>c_0$ there exists an $n_0$, such that for all $n>n_0$,''
etcetera.

\section{Summary of results}

We establish several properties
of random \maxtsat, random \maxksat, and random \maxcut,
focusing on 2-\sat.
This section summarizes our main results and indicates the
nature of the proofs;
further results and proofs are given in subsequent sections.

One of our goals is to establish the \maxtsat results
without depending on those for decision 2-\sat,
and in particular to work independently of~\cite{BBCKW01}
and reproduce its results.
We enjoy some success in this; the exceptions are our reliance
on~\cite{BBCKW01} for the upper bound in Theorem~\ref{MainLowDens}
(with an extraneous logarithmic factor arising in the translation),
and a more acute form of the same problem in the scaling window,
where we lack any corresponding bound for the $\l>1$
case of Theorem~\ref{MainScalingWindow}.

\begin{figure}[htbp]
\begin{center}
\psfrag{c}[cr]{density \fmath{c}}
\psfrag{0}[br]{$0$}
\psfrag{1}[br]{$1$}
\psfrag{3}[br]{$3/4$}
\psfrag{P}[cr]{\fmath{f(n,cn)/(cn)}}
\psfrag{n2}[bl]{$n \ra \infty$}
\psfrag{left}[bc][bc]{}
\psfrag{right}[bc][bc]{}
\includegraphics[width=4.5in]{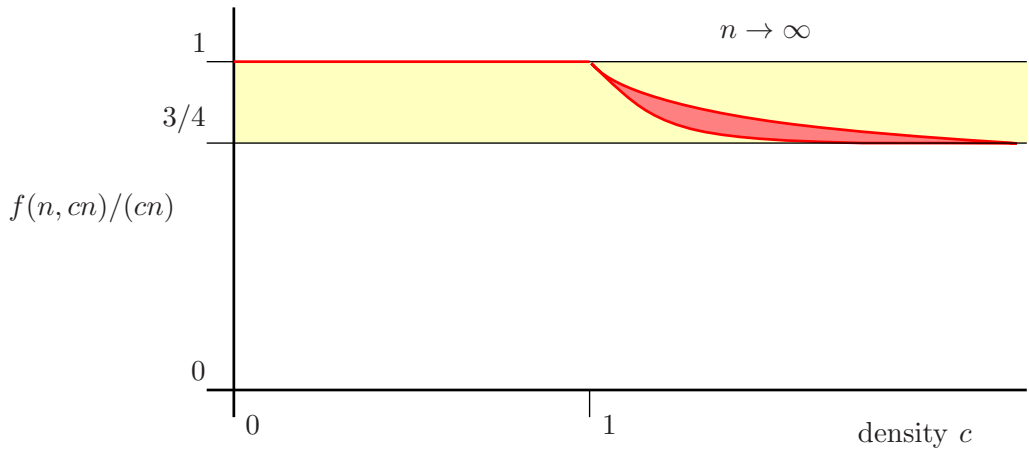}
\label{resultsfoil}
\end{center}
\caption{``Artist's rendition'' of the behavior of $f(n,cn)/(cn)$.}
\end{figure}

Figure~\ref{resultsfoil} show an ``artist's rendition'' of the
our results for 2-\sat.
For $c<1$, we expect to satisfy nearly all clauses,
while for $c \ra \infty$, we expect to satisfy only about 3/4ths of them.
The aysmptotic behavior for $c<1$ is understood;
so is that for $c$ large
(with a log-factor gap in the bounds on the second term);
and for $c=1 \pm \Theta(n^{-1/3})$
(with only a one-sided bound on the second term).
We now state these results more exactly;
we prove them in the next section.

For $c<1$ a random formula $F(n,cn)$ is satisfiable \whp,
so we would expect $\max F$ to be close to $cn$ in this case;
the following theorem shows this to be true.

\can{\subthreshold}{%
\begin{theorem}
\label{subthreshold}
For $c = 1-\e$, with any constant $\e>0$,
$\cn - f(n,\cn) = \Theta(1/(\e^3 n))$.
\end{theorem}
}
\uncan{\subthreshold}

The proof comes from
counting the expected number of the ``bicycles''
shown by \cite{CR92} to be necessary components
of an unsatisfiable formula.

For any $c$, $f(n,cn)\geq \frac34 cn$,
since a random assignment of the variables
satisfies each clause with probability $\frac34$.
The next theorem shows that neither this bound nor
the trivial upper bound $cn$ is tight,
although for large $c$, $\frac34 cn$ is close to correct.

\can{\mainhigh}{%
\begin{theorem}\label{MainHighDens}
For $c$ large,
$(\sqrt{c} \frac{\sqrt8 - \sqrt1}{3\sqrt{\pi}} -O(1))n
\leqa f(n,c n) - \frac34 cn \leqa
(\sqrt{c} \, \sqrt{3 \ln(2) / 8}) n $.
\end{theorem}
}
\uncan{\mainhigh}

\noindent
The values of
$\frac{\sqrt8 - \sqrt1}{3\sqrt{\pi}}$ and $\sqrt{3 \ln(2) / 8}$
are approximately $0.343859$ and $0.509833$, respectively.
The upper bound is proved by a simple first-moment argument,
and the lower bound by analyzing an algorithm;
both techniques are exactly those demonstrated
in~\cite[Lecture 6]{TenLectures}
to analyze the Gale-Berlekamp switching game.

Our next results relate to the low-density case,
when $c$ is above but close to the critical value~$1$.
How does $f(n,cn)$ depend on $c=1+\e$
for small $\e$? 

\can{\mainlow}{%
\begin{theorem}\label{MainLowDens}
For any fixed $\e>0$,
$(1+\ve- \e^3 /3)n
 \leqa f(n,(1+\e)n)$;
also,
there exist absolute constants $\a_0$ and $\e_0$,
such that for any fixed $0< \e < \e_0$,
$f(n,(1+\e)n) \leqa
 (1+\e -{1\over 3}{\a_0 \e^3 / \ln(1/\e)}) n$.
\end{theorem}
}
\uncan{\mainlow}

That is, a constant fraction of the clauses must remain unsatisfied,
but this fraction --- $\e^3/3$ at most --- is surprisingly small.
The lower bound is proved by using the ``differential equation method''
(see for example~\cite{worm})
to exactly analyze
a version of the unit-clause heuristic.
The upper bound's proof is a simple first-moment argument;
however, for the probability that a sub-formula with density $>1$
is satisfiable,
it requires the
exponentially small bound
given by \boll et al.~\cite{BBCKW01}
(see Theorem~\ref{littleboll} below).
It is likely that,
by replacing our use of \cite{BBCKW01}
with structural properties of the kernel of a sparse random graph,
the upper bound's $\e^3/\ln(1/\e)$ can be replaced by $\e^3$
to match the lower bound up to constants%
~\cite{JS02a}.

The major significance of \cite{BBCKW01}
was to determine the
``scaling window'' for random 2-\sat.
Without using their result,
we prove an analogous result for \maxtsat,
and incidentally reproduce most parts of their 2-\sat result.

\can{\mainscaling}{%
\begin{theorem}\label{MainScalingWindow}
Letting $c = 1+\epsilon = 1 + \lambda n^{-1/3}$,
we have
\begin{align*}
\cn  - f(n,\cn) &=
\begin{cases}
O(\lambda^3) & \text{if $\lambda > 1$;}
\\
\Theta(1)    & \text{if $-1 \leq \lambda \leq 1$;}
\\
\Theta(|\lambda|^{-3}) & \text{if $\lambda < -1$.}
\end{cases}
\end{align*}
Furthermore, for $\l>1$, for some positive absolute constant, and any $\b>0$,
\begin{align*}
 \Pr\left((\cn-f(n,\cn)) > \const \b \l^3\right)
  & \leq \exp(-\b \l^3) .
\end{align*}
Also,
\begin{align*}
\Pr(F(n,cn) \text{ is satisfiable}) &=
\begin{cases}
\exp(-O(\lambda^3)) & \text{if $\lambda > 1$;}
\\
\Theta(1) & \text{if $-1 \leq \lambda < 1$;}
\\
1-\Theta(|\lambda|^{-3}) & \text{if $\lambda < -1$.}
\end{cases}
\end{align*}
\end{theorem}
}

\uncan{\mainscaling}

In particular, in the scaling window $c = 1 \pm \lambda n^{-1/3}$,
a random formula is satisfiable with probability
which is bounded away from
0 and~1 (the exact bounds depending on $\lambda$),
and it can be made satisfiable by removing a constant-order
number of clauses (the constant depending on $\lambda$).

In section~\ref{sec:ksat},
for \maxksat, we derive analogous results only for $c$ large,
reflecting the general state of ignorance
regarding the $k$-\sat phase transition.
(For some results on scaling windows for $k$-\sat see~\cite{WILSON}.)
Still more generally, Theorem~\ref{csp} describes the high-density
case for any \maxcsp.
More interestingly,
for random \maxksat (including $k=2$)
we observe that $\max F$
is concentrated about its expectation $f(n,cn)$
(as previously remarked in~\cite{BFU93})
and that $f(n,cn)/(cn)$ is monotone non-increasing in~$c$.
Were $f(n,cn)/(cn)$ also monotone in $n$, an important property
analogous to the satisfiability conjecture would follow;
we present this as a conjecture for general \maxcsp{s}.

In section~\ref{sec:online} we consider online versions of \maxtsat,
for one of which we prove that a natural greedy algorithm is optimal.

Results for the \maxcut problem for sparse random graphs,
which is closely analogous to random \maxtsat,
are presented in section~\ref{sec:cut}.

\section{Random MAX 2-SAT}

\subsection{Sub-critical MAX 2-SAT}
One of the most basic facts concerning \maxtsat is that for
constants $c<1$, the expected number of clauses unsatisfied
is~$o(1)$.
This is refined by Theorem~\ref{subthreshold}, which shows the
number to be~$\Theta(1/(\e^3 n))$.
We now prove the theorem.

\thmpf{subthreshold}{Proof.}
We write the proof in the \sat equivalent of the ``\gnp'' model,
because the expressions for the probability of a clause's presence are
cleaner in this model, but adaptation to the \gnm model is immediate.

\begin{figure}[tbp]
\begin{center}
\psfrag{x1}[tl][tl]{\fmath{u=w_0}}
\psfrag{x1b}[br][br]{\fmath{w_i}}
\psfrag{xk}[tr][tr]{\fmath{w_j}}
\psfrag{xkb}[br][br]{\fmath{v=w_{k+1}}}
\includegraphics[width=3in]{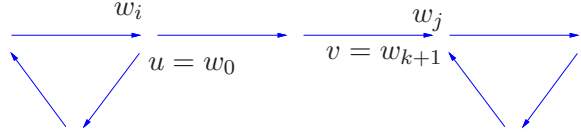}
\end{center}
\label{bicycle}
\caption{Sequence of clause-derived implications for a bicycle.
Start the walk from $u$,
proceed clockwise to $w_i$ (which equals either $u$ or $\bar{u}$),
continue right to $w_j$,
and again go clockwise to terminate at $v$
(which equals either $w_j$ or $\bar{w}_j$).}
\end{figure}

A $k$-\emph{bicycle}
(see Figure~\ref{bicycle})
is a sequence of clauses
$\set{\bu,w_1},$ $\set{\bw_1,w_2},$ $\ldots,$ $\set{\bw_k,v}$ where
literals $w_1,w_2,\ldots,w_k$ are distinct \emph{as variables}
(none is the same as nor the complement of another)
and $u\in \set{w_i,\bw_i}$,
$v\in\set{w_j,\bw_j}$
for some $1\leq i,j\leq k$.
(Think of it as a ``walk'' in which the first and last variables
are also both visited \textit{en route}.)
Because satisfying a clause $\set{\bu,v}$ means that if $u$ is true then
$v$ must be true, such a clause yields an implication
$u \ra v$ (and a complementary implication $\bv \ra \bu$);
Figure~\ref{bicycle} represents such a sequence of implications
for a bicycle.
Chv\'atal and Reed~\cite{CR92} argue that if a formula is infeasible then it
contains a bicycle.
Thus if we delete an edge from every bicycle,
the remaining subformula is satisfiable.

The number of potential $k$-bicycles, whether or not present in a
given formula~$F$, is at most $k^2 (2n)^{k}$. The probability
that all $k+1$ clauses of a given bicycle are present in a random formula $F$ is at most
$[(cn) / (2^2 \binom{n}{2})]^{k+1} = [c/(2(n-1))]^{k+1}$, so the
expected number of $k$-bicycles is $\leqa k^2 c^{k+1} / (2n)$.
If we delete one edge in every bicycle, we obtain a satisfiable formula.
For any fixed $c<1$, $\sum_{k=1}^n k^2 c^{k+1} / (2n) = \Theta(1/n)$.
Thus, the expected number of edges we need to delete is at most
$O(1/n)$ and $f(n,\cn) \geq \cn - O(1/n)$.

To obtain the lower bound we show that with probability at least
$\Theta(1/(\e^3 n))$
the formula $F$ is not satisfiable. This clearly implies an upper bound
$f(n,\cn) \leq \cn - \Theta(1/(\e^3 n))$. To this goal
we employ the second moment method.

For simplicity here, we will restrict
ourselves to $3$-bicycles,
which will only establish ``$\Theta_{\e}(1/n)$'',
that is, something of order $\Theta(1/n)$ but with hidden constants
that may depend on~$\e$.
The full proof is the same but using bicycles of lengths
up to $1/\e$, not just length~3,
and parallels the proof of Theorem~\ref{MainScalingWindow},
case $\l \leq -1$.
(In fact, taking $\l = \l(n) = \e n^{1/3}$ there
establishes the current theorem completely.)

Consider $4$-tuples of  clauses of the form
$\set{\bu_1,u_2},\set{\bu_2,\bu_1},\set{u_1,u_3},\set{\bu_3,u_1}$, where
$u_1,u_2,u_3$ are arbitrary variables. One observes that this sequence
of clauses is a $3$-bicycle, and, moreover, its presence in the random
formula~$F$ implies non-satisfiability. We now show, using second moment
method, that the number $B_3$ of such bicycles is at least one with probability at least $\Omega(1/n)$.
We have $\E[B^2_3]=\sum \pr(X\in F,X'\in F)$,
where the sum runs of over the pairs of $3$-bicycles $X,X'$ of the form above,
and $X\in F$ means all the clauses of $X$ are present in $F$. We decompose the
sum into three parts: the sum over pairs $X,X'$ with $X=X'$, the sum over pairs
that do not have common clauses and the rest. It is easy to see that the first
sum is simply $\E[B_3]$ which is $\Theta(1/n)$, by the argument for upper bound.
To analyze the second sum note that for each fixed pair $X,X'$ with no common clauses,
we have $\pr(X,X'\in F)=\pr(X\in F)\pr(X'\in F)$, when replacement of clauses is allowed.
(When replacement is not allowed the reader can check that the difference between the
left and the right-hand sides is very small, and the rest of the argument goes through).
Then, this sum is smaller than  $\sum_{X,X'}\pr(X\in F)\pr(X'\in F)=(\E[B_3])^2$,
where the sum now runs over all the pairs $X,X'$. For the third sum we have two
cases. First case is pairs $X,X'$ defined on the same set of variables. For
example $X=\set{\bu_1,u_2},\set{\bu_2,\bu_1},\set{u_1,u_3},\set{\bu_3,u_1}$
and $X'=\set{\bu_1,u_2},\set{u_2,u_1},\set{\bu_2,u_3},\set{\bu_3,\bu_2}$, share
one clause $\set{\bu_1,u_2}$ and are defined over the same set of variables.
There are $O(n^3)$ choices for the variables $u_1,u_2,u_3$ in these pairs.
But since $X\neq X'$ then there are altogether at least five clauses in $X$
and $X'$ together. For a given pair, the probability that all these clauses
are present in $F$ is $O(1/n^5)$. Then the expected number of such pairs $X,X'\in F$
is $O(1/n^2)=o(1/n)$.

The second case is pairs $X,X'$ defined over different set of variables.
Since they share a clause then the pair is defined on exactly four variables.
But then there are at least six clauses in this pair. We obtain that the
expected number of such pairs $X,X'$ which belong to $F$ is at most $O(n^4)O(1/n^6)=O(1/n^2)=o(1/n)$.

We conclude that $\E[B_3^2]=\E[B_3]+o(1/n)=\Theta(1/n)+o(1/n)$. We now use the bound
$\pr(Z\geq 1)\geq (\E[Z])^2/\E[Z^2]$, which holds for any non-negative integer random variable $Z$.
Applying this bound to $B_3$ we obtain $\pr(B_3\geq 1)\geq (\E[B_3])^2/\E[B_3^2]\geq
\Theta(1/n^2)/(\Theta(1/n)+o(1/n))=\Theta(1/n)$. This completes the proof.
\qed

It is worth pointing out the following simple fact,
upon which we will shortly improve.

\begin{remark}
\label{basics}

For $c>1$, $f(n,cn) \geqa n (\frac34 c + \frac14)$.
\end{remark}

\begin{proof}
It suffices to show that for any $\e>0$,
for all $n$ sufficiently large, $f(n,cn) \geq (\frac34 c+\frac14-\e)n$.
Select the first $(1-\e)n$ clauses,
and let $\Xvec$ be a best assignment for it.
By Theorem~\ref{subthreshold},
$\Xvec$ satisfies an expected $(1-\e)n -o(1)$
of these first clauses. Also,
an expected $3/4$ths of the remaining $(c-1+\e)n$ clauses are satisfied,
yielding the claim.
\end{proof}

\subsection{High-density random MAX 2-SAT}

While it is well known that for $c>1$, $F(n,cn)$ is a.s.\ unsatisfiable,
is it possible that even for $c$ large, \emph{almost} all clauses are
satisfiable?
Theorem~\ref{MainHighDens} rules this out by
showing that a constant fraction
of clauses must go unsatisfied;
up to a constant, it also provides a matching lower bound.

\thmpf{MainHighDens}{Proof of the upper bound.}
The proof is by the first-moment method.
If $\max F > (1-r)cn$ then there is a satisfying assignment of a
subformula $F'$ which omits $rcn$ or fewer clauses,
and where (taking $F'$ to be maximal)
all the omitted clauses are unsatisfied.
Any fixed assignment satisfies each (random) clause of $F'$ w.p.~$3/4$
and unsatisfies each omitted clause w.p.~$1/4$,
so by linearity of expectations,
the probability that there exists such an $F'$ is

\begin{align*}
 P = \pr(\exists \text{ satisfiable } F')
 & \leq
 2^n \sum_{k=0}^{r c n} \binom{cn}{k} (\frac34)^{cn-k} (\frac14)^k .
\end{align*}
For $r<\frac14$ the sum is dominated by the last term.
From Stirling's formula $n! \asymp \sqrt{2 \pi n} \; (n/e)^n$,
\begin{align}
\label{binom}
\binom{cn}{r c n}
 & \asymp
 1/\sqrt{2 \pi r(1-r)cn} \quad (r^{-r} (1-r)^{-(1-r)})^{cn} .
\end{align}

Substituting \eqref{binom} into the previous expression,
\begin{align*}
 P & \leqa
 1/\sqrt{2 \pi r(1-r)cn}2^n \quad(r^{-r} (1-r)^{-(1-r)} (3/4)^{1-r} (1/4)^r)^{cn} .
\end{align*}
Substituting $r=1/4-\e$,
\begin{align*}
 \frac{1}{cn} \ln P
 & \leqa
 \ln(2)/c - (8/3) \e^2 + O(\e^3) ,
\end{align*}
so that for $\e > \sqrt{(3/8) \ln 2 / c}$,
as $n \ra \infty$, $P \ra 0$.
The conclusion follows.

\pf{Proof of the lower bound.}  The proof is algorithmic.
When variables $X_1,\ldots,X_{k}$ have been set, define the
reduced formula $F_{k}$ in which any clause containing a True literal
is removed and ``scored'',
and False literals are removed from the remaining clauses.
(Clauses with 0 variables remaining are permanently unsatisfied.)
Define a potential function $q(F_{k})$ to be the number of clauses
already satisfied,
plus $3/4$ the number of 2-variable clauses (``2-clauses''),
plus $1/2$ the number of 1-variable clauses (``unit clauses'').
Note that randomly assigning the remaining variables satisfies
an expected total number of clauses precisely $q(F_{k})$, so
$q$ is a lower bound on the number of clauses satisfiable.

After variables $X_1, \ldots, X_{k-1}$ have been set to define $F_{k-1}$,
our algorithm sets $X_k$ in whichever of the two ways gives an $F_k$ with
larger value $q(F_k)$.  (Ties may be broken arbitrarily.)
In $F_{k-1}$, let the number of appearances of $X_k$ and $\Xbar_k$
in unit clauses be denoted by $A_1$ and $\Abar_1$,
and their number of appearances
in 2-clauses by $A_2$ and $\Abar_2$.
If $X_k$ is set to True, then
\begin{align*}
q(F_k) - q(F_{k-1})
 &= \Delta_k
 \eqdef
 \frac12 (A_1 - \Abar_1) + \frac14(A_2 - \Abar_2),
\end{align*}
and if $X_k$ is set False, then $q(F_k) - q(F_{k-1}) = -\Delta_k$.
Note that
$q(F_k) = q(F_{k-1}) + |\Delta_k|$
is a lower bound on the number of satisfiable clauses,
and $q(F_0) = \frac34 cn$.

With $k-1$ variables already set, $F_{k-1}$ a.s.\ has
$(\frac12 \pm O(1/\sqrt{c})) 2 \frac{k-1}{n} \frac{n-k+1}{n} \cdot cn$
unit clauses, and
$\left(\frac{n-k+1}{n}\right)^2 \cdot cn$ 2-clauses,
on the remaining variables.
(The reason for $\frac12 \pm O(1/\sqrt{c})$ instead of $\frac12$ is that
we set the previous variables in a biased manner.)
Also, conditioned on the number of clauses,
$F_{k-1}$ is a uniformly random formula (each ``slot'' being equally
likely to be filled by any of the remaining literals).
For $n$ large,
$A_1$ and $\Abar_1$ are approximated by independent Poisson
random variables with parameter $(\frac12 \pm O(1/\sqrt{c})) \frac{k-1}{n} c$,
and $A_2$ and $\Abar_2$ by Poissons with parameter $ \frac{n-k+1}{n}c$.
By assumption, $c$ is large,
so each of these distributions is approximately Gaussian,
and their sum $\Delta_k$ is also approximately Gaussian,
with mean~0 (by symmetry)
and variance
\begin{align*}
\sigma_k^2 &=
2 \cdot (\frac{1}{2})^2 \cdot \Var(A_1) +
2 \cdot (\frac{1}{4})^2 \cdot \Var(A_2)
\\ &
= c \left( (\frac14 \pm O(1/\sqrt{c})) \frac{k-1}{n}
           + \frac18  \frac{n-k+1}{n} \right) .
\end{align*}
For $Z \distr N(0,1)$, it is well known that $\E|Z| = \sqrt{2/\pi}$;
thus $\E|\Delta_k| = \sqrt{2/\pi} \: \sigma_k=
\sqrt{2/\pi}\sqrt{c(\frac14 \frac{k-1}{n} + \frac18  \frac{n-k+1}{n})} \pm O(1)$.
Finally,
\begin{align*}
\E(q(F_n)) & \geq \frac34 c n + \sum_{k=0}^{n-1} \E(|\Delta_k|)
\\ &
 \approx
 \frac{3}{4}cn+ \int_{0}^{n} \E(|\Delta_k|) dk
\\ &
\geqa
 \frac{3}{4}cn+
 \left( \sqrt{c} \frac{\sqrt{8}-\sqrt{1}}{3 \sqrt{\pi}} - O(1) \right) n .
\end{align*}
\qed

We remark that
in the preceding proof, $X_k$ was set \true or \false so as to maximize
half the number of satisfied unit clauses plus
a quarter the number of satisfied 2-clauses.
This is reminiscent of the ``policies'' in~\cite{AchSo00}.
There, the goal was to satisfy as a dense a 3-\sat formula as possible;
unit clauses always had to be satisfied,
and variables were set so as to maximize a linear combination of the
number of satisfied 2-clauses and 3-clauses.
In~\cite{AchSo00}, the linear combination which was optimal for the
purpose changed during the course of the algorithm;
the determination of the optimal combinations, and the proof of
optimality, was a main result of the paper.
In the present case, though, it is evident that the ratio 1:2 is optimal:
for $c$ large, the potential function $q$ predicts the expected number of
clauses satisfiable almost exactly.
The difference can be ascribed to the fact that
here $c$ is ``large'', and in~\cite{AchSo00} the corresponding parameter
(the initial 3-clause density) was fixed
(relevant values were in the range of 3.145 to 3.26).
Were we to try to tune the \maxtsat algorithm above for small values of $c$,
more complex methods like those of~\cite{AchSo00} would presumably be needed.

\subsection{Low-density random MAX 2-SAT}

For low-density formulas, with $c = 1+\e$ and $\e>0$ a small constant,
the bounds of Theorem~\ref{MainHighDens} are inapplicable.
It is still true (from Remark~\ref{basics})
that we expect to satisfy at least $(1+ \frac34 \e)n$ clauses,
but it is not obvious whether
the best answer is this, or close to the full number of clauses
$(1+\e)n$, or something in between. In this section we prove
Theorem~\ref{MainLowDens}  which shows that $(1+\e)n - f(n,cn)$,
the number of clauses we must dissatisfy,
lies between $\Theta(\e^3 n/\ln(1/\e))$ and $\Theta(\e^3 n)$.
That is, a linear fraction of clauses must be rejected,
but this fraction, at most $\Theta(\e^3)$, is surprisingly small.
We will employ the following theorem of \boll et al.~\cite{BBCKW01}
on random 2-\sat.

\begin{theorem}\textnormal{(\cite{BBCKW01}, Corollary 1.5)}
\label{littleboll}
There exist positive constants $\a_0$ and $\e_0$
such that for any $0 < \e < \e_0$
and sufficiently large $n$,
$\pr[ F(n,(1+\e)n) \text{ is satisfiable}] \leq
 \exp(-\alpha_0 \e^3 n)$.
\end{theorem}
\nobreak \noindent%
(Here, $\alpha_0$ is the $\liminf$ of the constant implicit in $\Theta$ in
the theorem in~\cite{BBCKW01}.)
The $\exp(-\Theta(\e^3 n))$ probability of satisfiability in
random 2-\sat translates into an expected $O(\e^3 n / ln(1/\e))$
unsatisfied clauses in random \maxtsat.

\par
\vspace{.1in}
\thmpf{MainLowDens}{Proof of the upper bound.}
The proof is by the first-moment method.
Let $c=1+\e$.
Let $F'$ range over subformulas of $F$ which omit $r c n$ or fewer clauses.
Specifying $r < 1/4$,
the conditions of Theorem~\ref{littleboll} apply, so
\begin{align}
\label{logbound}
 P = \pr(\exists \text{~maximally~ satisfiable } F')
 & \leq
 \sum_{k=0}^{r c n} \binom{cn}{k} ({1\over 4})^{r c n}e^{-\alpha_0(\e -\frac{k}{n})^3n} ;
\end{align}
as $r<1/4$, the sum is dominated by the last term.
Using~\eqref{binom} to approximate $\binom{cn}{crn}$,
\begin{align*}
\frac{1}{c n} \ln P
 & \leqa
-\logis{r} - \logis{(1-r)} - \a_0 (\e-c r)^3/c - r\ln(4).
\\ \intertext{First observe that as $\e \ra 0$, for any $r = o(\e)$,
this is}
 & = -r \ln r (1+o(1)) - \a_0 \e^3 (1+o(1))  - r\ln(4) .
\\ \intertext{For any constant $b<1/3$,
if $r = b \a_0 \e^3 / \ln(1/\e)$,
this is}
 & =
 3b \a_0 \e^3 (1+o(1)) - \a_0 \e^3 (1+o(1))<0.
\end{align*}
That is, it is unlikely that asymptotically fewer than
$(1/3)\a_0 \e^3 / \ln(1/\e)$ clauses can go unsatisfied.

\pf{Proof of the lower bound.} The proof is algorithmic,
and of the sort familiar from~\cite{AchSo00} and previous works.
It analyzes a version of the ``unit-clause'' heuristic.
Initially, ``seed'' the algorithm by randomly deleting a variable
from each of, say, $n^{1/10}$ random 2-clauses
to convert them to unit clauses.
While $F$ has any unit clauses, select one at random and set its variable
to satisfy the clause.
Continue until no unit clauses remain.
The analysis consists of counting the clauses unsatisfied in these steps,
and justifying the assertion that when there are no more unit clauses,
$o(1)$ further clauses need be unsatisfied.

When $k$ variables have been set, let the number of 2-clauses
be denoted $m_2(k)$, the number of unit clauses $m_1(k)$,
and the number of unset variables $m(k) = n-k$.
In one step, the changes in these quantities are $\D m = -1$,
$\E(\D m_2) = -\frac{2}{m} m_2$, and
$\E(\D m_1) = -1 - \frac{1}{m} m_1 + \frac{1}{m} m_2 $
(assuming that $m_1 >0$ before the step).
Over a large number of steps, the net changes will be a.s.\ a.e.\
equal to the expectations.
Renormalizing with $\rho = m/n$, $\rho_1 = m_1/n$, and $\rho_2 = m_2/n$,
the differential equation method (see for example~\cite{AchSo00,worm})
asserts that $(\rho_1,\rho_2)$
a.s.\ a.e.\ obey the differential equations
\begin{align*}
 d \rho_2 / d\rho &= \frac{2 \rho_2}{\rho}
 &
 d \rho_1 / d\rho &= 1 + \frac{\rho_1}{\rho} - \frac{\rho_2}{\rho} .
\\
\intertext{%
With boundary conditions that for $\rho=1$ (i.e., initially),
$\rho_2 = c$ and $\rho_1 = 0$,
the unique solution is
}
 \rho_2 &= c \rho^2
 &
 \rho_1 &= c \rho - c \rho^2 + \rho \ln \rho .
\end{align*}
This results in $\rho_1=0$ at two times: initially, when $\rho=1$,
and also for $\rho = \rhostar$ satisfying
\begin{align}
 c &= \ln(\rhostar) / (\rhostar-1) .
 \label{rhostar}
\end{align}

While $\rho>\rhostar$, the only clauses ever unsatisfied are unit
clauses which contain the negation of the variable being set,
and the expected number of such rejected clauses per step is
$\frac{1}{2m} m_1 = \frac{\rho_1}{2\rho}$.
Integrating over the period $\rhostar$ to~1,
\begin{align}
\int_{\rhostar}^1 \frac{\rho_1}{2 \rho} d\rho
  &=
\frac12 \int_{\rhostar}^1 \left( c - c \rho + \ln \rho \right) d\rho
\notag
\\ &=
\left.
  \frac12
  \left( c \rho - c \rho^2 / 2  + \rho \ln \rho  - \rho \right)
\right|_{\rhostar}^{1}
\notag
\\
  \intertext{which, substituting for $c$ from~\eqref{rhostar}}
  &=
\frac12 (\rhostar-1) - \frac14(\rhostar+1) \ln \rhostar .
\label{int1}
\end{align}

So from $\rho=1$ to $\rho=\rhostar$, the number of clauses dissatisfied
by the algorithm is a.s.\ a.e.\ $n$ times expression~\eqref{int1}.
After this time, the remaining (uniformly random) 2-\sat formula has density
$\rho_2(\rhostar) \; / \; \rhostar
  = c \rhostar
  = \ln(\rhostar)\rhostar/ (\rhostar-1)
  < 1
$
since $\ln(\rhostar)<\rhostar-1$ and $\rhostar<1$,
and thus (by Theorem~\ref{basics}) contributes $o(1)$ to the expected
number of unsatisfied clauses.
In short, the algorithm a.s.\ fails to satisfy a.e.\
$(\frac12 (\rhostar-1) - \frac14(\rhostar+1) \ln \rhostar) n$ clauses.
For $\rhostar$ (asymptotically) close to~1,
the number of dissatisfied clauses is $\asymp n(1- \rhostar)^3/24$.
In particular,
with $\e>0$ asymptotically small and $c = 1+\e$,
$\rhostar \asymp 1-2\e$,
and the number of dissatisfied clauses is $\asymp n\e^3/3$.
 \qed

Two remarks.  First,
in addition to the asymptote, the proof gives
a precise parametric relationship (as functions of $\rhostar$)
between the clause density $c$
(given by~\eqref{rhostar})
and the rejected-clause density
(given by~\eqref{int1}).
Solving numerically,
for $c=1.5$ we find
rejected-clause density $\approx 0.0183275$,
and for $c=2$
--- where naively the rejected-clause density would be
$\frac14 c = 0.5$ --- we achieve
rejected-clause density $\approx 0.0809517$.

Second, with the solution in hand, the asymptotic behavior
is easy to see without the need for differential equations.
This alternate proof is not fully rigorous, but is more intuitive and more robust;
it is the basis of the analysis within the scaling window
(see Theorem~\ref{MainScalingWindow}).

\thmpf{MainLowDens}{Alternate proof of lower bound.}
Consider what happens when $m = (1-\d)n$ variables remain unset.
The number of 2-clauses is a.s.\
$m_2 \asymp (1-\d)^2 (1+\e) n \asymp (1+\e-2\d)n$.
The expected increase in the number of unit clauses is then
$ \E(\D m_1) =  -1 - m_1/m + m_2/m \geq -1+ m_2/m$
(and the neglected $m_1/m$ is not only conservative, but will
also prove to be insignificantly small).
Thus,
$\E(\D m_1) \geq -1+[(1+\e-2\d)n] / [(1-\d)n]
 \asymp \e-\d$.
At $\d=0$, the number of unit clauses increases by
$\e$ per step, this increase linearly falls to 0 per step by $\d=\e$,
and further to $-\e$ by $\d=2\e$:
the expected number of unit clauses is bounded by an inverted parabola,
with base $2\e n$ and height $\frac12 \e^2 n$.
At each step about $1/(2n)$th of the unit clauses get dissatisfied.
The area under the parabola, times this $1/(2n)$ factor, is
$\frac23 \cdot \text{base} \cdot \text{height} \cdot 1/(2n)
 = \frac13 \e^3 n$.
\qed

\section{The MAX 2-SAT scaling window}
For random \maxtsat, we have seen that
for fixed $c<1$, $\cn - f(n,\cn) = \Theta(1/n)$,
and for $c>1$, $cn-f(n,cn) = \Theta(n)$.
That is, random \maxtsat experiences a phase transition around $c=1$.
It is natural to ask about the scaling window around the
critical threshold: What is the interval around $c=1$ within which
$\cn - f(n,\cn) = \Theta(1)$?
Theorem~\ref{MainScalingWindow} shows that the scaling window is
$c = 1 \pm \Theta(n^{-1/3})$.
The corresponding question for random 2-\sat is the range in which
$\pr(F(n,\cn) \text{ is satisfiable}) = \Theta(1)$.
This was shown by~\cite{BBCKW01} to be $c = 1 \pm \Theta(n^{-1/3})$
with their result reproduced as Theorem~\ref{boll} here.

\begin{theorem}[\textbf{\boll et al, \cite{BBCKW01}}] \label{boll}
Let $F(n,cn)$ be a random 2-\sat formula, with
$c=1+\lambda_n n^{-{1 / 3}}$.
There are absolute constants
$0<\e_0<1$,
$0<\lambda_0<\infty$,
such that the probability $F$ is satisfiable is:
$1-\Theta({1/ |\lambda_n|^3})$,
when $-\e_0 n^{1/ 3}\leq \lambda_n\leq -\lambda_0$;
$\Theta(1)$, when $-\lambda_0\leq \lambda_n\leq \lambda_0$;
and
$e^{-\Theta(\lambda_n^3)}$,
when $\lambda_0\leq \lambda_n\leq \e_0 n^{1/ 3}$.
\end{theorem}

That the two scaling windows are the same is no coincidence,
and in fact Theorem~\ref{MainScalingWindow}
reestablishes much of Theorem~\ref{boll} independently.

\thmpf{MainScalingWindow}{Proof.} Note that, provided we prove the bounds
for the cases $\lambda\leq -1$ and $\lambda\geq 1$,
the bound for the case $|\lambda|<1$ follows immediately,
since we obtain that the probability of satisfiability is at least
$\exp(-O(\lambda^3))\geq \exp(-O(1))$
and at most $1-\Theta(1/|\lambda|^3)\leq 1-\Theta(1)$,
where in both cases $|\lambda|<1$ was used.
The more interesting cases $|\lambda|\geq 1$
are considered in two subsections below.

\subsection{Case $c=1+\lambda n^{-1/3}$, $\lambda\leq -1$}
For convenience we write $c=1-\lambda n^{-1/3}$ and $\lambda\geq 1$.
The proof for this case is very similar to that
of Theorem~\ref{subthreshold} and uses the notion of bicycles.
(As in the earlier case, we work in the equivalent of the \gnp
model for notational convenience, with the understanding that
the proof works equally well in the \gnm model.)
As before, the number of clauses that must be dissatisfied
is bounded by the number of bicycles.
The expected number of $k$-bicycles is at most $k^2c^{k+1}/(2n)=
k^2(1-\lambda n^{-1/3})^{k+1}/(2n)$.
Using the formula
$\sum_{k\geq 1}k^2\rho^k={\rho(1+\rho)\over (1-\rho)^3}$
which for $\rho \asymp 1$ is $\asymp 2/(1-\rho)^{3}$,
we have
\begin{align}
\sum_{1\leq k<\infty} k^2 (1-\lambda n^{-1/3})^{k+1}/(2n)
 & \asymp 2/\lambda^3 .
 \label{bound1}
\end{align}
Therefore $\lfloor cn \rfloor - f(n,\lfloor cn \rfloor)=O(1/\lambda^3)$.
Using Markov's inequality we also
obtain that the probability that the formula is unsatisfiable
is at most the expected number of bicycles,
that is, at most $O(1/\lambda^3)$.

We now obtain a matching lower bound. Consider only ``bad'' bicycles,
in which $u=\bw_i$, $v=\bw_j$, and $i<j$.
Note that no bad bicycle is completely satisfiable,
since the first ``wheel''
$u \ra \cdots \ra \w_i = \bu$
requires $u = \false$
and thus $w_i = \true$;
whereupon the path
(technically called the ``top tube'' of a bicycle)
$w_i \ra \cdots \ra \w_j$ implies $w_j = \true$;
and the second wheel
$w_j \ra \cdots \ra v = \bw_j$ provides a contradiction.
Note that about $1/8$th of the potential bicycles are bad.

Let $B_k$ denote the number of the bad $k$-bicycles.
Since
\begin{align}
\E(\text{\#unsatisfiable clauses})
 &\geq \Pr(\text{$F$ unsatisfiable})
 \\ & \geq \Pr(\sum_{k \leq K} B_k \geq 1) ,
 \notag
\\ \intertext{it suffices to prove that this is}
 & = \Omega(1/\lambda^3) ;
 \label{desired1}
\end{align}
we will show this for $K= (1/\lambda)n^{1\over 3}$.
Repeating the argument for~\eqref{bound1}, we obtain that
\begin{align*}
\E[\sum_{k\leq K}B_k]
 & \geqa (2/(8e)) / \l^3,
\end{align*}
the $1/(8e)$ coming from the series' truncation at $K$ and
the use of only bad bicycles.
To obtain~\eqref{desired1} it suffices prove that
\begin{align}
\E[(\sum_{k \leq K} B_k)^2]
 & =
 (1+O(1)) \cdot \E[\sum_{k \leq K} B_k] ,
\label{SecondMoment}
\end{align}
for then
\[
\pr(\sum_k B_k\geq 1)
\geq
{(\E[\sum_k B_k])^2\over \E[(\sum_k B_k)^2]}
=
{(\E[\sum_k B_k])^2\over \E[\sum_k B_k] (1+O(1))}
=\Omega(1/\lambda^3) .
\]

We will prove~\eqref{SecondMoment} with $O(1/\l^3)$ filling
in for $O(1)$ (recall that $\l \geq 1$).
Consider pairs of $k,k'$-bicycles $X,X'$
with $k,k' \leq K$.
It suffices to show that for every~$X$,
\begin{align}
\label{conditioning}
\sum_{X' \neq X} \pr(X'\clausein F|X\clausein F) = O(1/\lambda^3) ,
\end{align}
because then
\begin{align*}
\E[(\sum_k B_k)^2]
 &= \sum_{X,X'}\pr(X,X' \clausein F)
 \\ & =
 \sum_X \Pr(X \clausein F) \;
  [1 + \sum_{X' \neq X} \Pr(X' \clausein F \given X \clausein F)]
 \\ & \leq
 \E[\sum_k B_k] (1 + O(1/\l^3)) .
\end{align*}

Establishing \eqref{conditioning} is the nub of the proof.
First, observe that for any bicycle $X'$ sharing no literals with~$X$,
$\Pr(X' \clausein F \given X \clausein F) \leq \Pr(X' \clausein F)$,
and so such bicycles $X'$ contribute $\leq \E \sum_k B_k = O(1/\l^3)$
to the sum.

Given a bicycle
$X'=\set{u,w_1},\set{\bw_1,w_2},\ldots,\set{\bw_k,v}$, a sequence
of literals $w_i,w_{i+1},\ldots,w_j$ from $X'$ is defined to be
a type~I excursion if literals $w_i,w_j$ belong to $X$ but literals
$w_{i+1},\ldots,w_{j-1}$ do not.
(If $j=i+1$, a sequence $w_i,w_{i+1}$ is a type~I
excursion if the corresponding clause $(\bw_i,w_{i+1})\in$ $X'$
does not belong to $X$.)
A sequence of literals $u',w_1,\ldots,w_j$
in $X'$ is defined to be a type~II excursion
if the literal $w_j$ belongs to $X$, but
$u,w_1,\ldots,w_{j-1}$ do not.
Similarly, a sequence $w_j,w_{j+1},\ldots,v'$ in $X'$
is defined to be a type~III excursion.

Bicycles $X'$ which are neither equal to $X$ nor disjoint from $X$
must have at least one excursion
(and at most one each of excursions of type~II and~III).
It suffices to establish \eqref{conditioning} for such
bicycles~$X'$.
We will just show that the expected number of bicycles $X'$
with one type~II excursion, no type~III excursion,
and any number $r \geq 0$ of type~I excursions,
is $O(1/\l^3)$;
the other three cases
(classified by the number of type~II and~III excursions)
follow similarly.

Since a collection of excursions uniquely defines~$X'$,
it is enough to count such collections.
Let the lengths of the type~I excursions be
$m_1,m_2,\ldots,m_r \geq 2$ and
that of the type~II excursion $m_{II}$,
where the length is defined by the number of literals.

For each type~I excursion there are two endpoints (literals)
which belong to $X$.
Since the size of $X$ is $\leq K= (1/\lambda)n^{1\over 3}$,
there are $\leq K^{2r} = (1/\lambda^{2r})n^{2r\over 3}$
choices for all the end points.
The $i$th type~I excursion contains $m_i-2$ literals not from $X$,
so
there are at most $(2n)^{m_i-2}$ ways of selecting them.
The excursion contains $m_i-1$ clauses,
all not from~$X$, so the probability they are all present in $F$
is $(1-\lambda n^{-1/3})^{m_i-1}/(2n)^{m_i-1}$.

Similarly, for the type~II excursion,
there are at most $K$ choices for the endpoint literal $w_{j-1}$,
which belongs to~$X$,
and at most $(2n)^{m_{II}-2}$ choices for
other literals $u',w_1,\ldots,w_{j-2}$.
The excursion contains $m_{II}-1$ clauses,
all not from~$X$, so the probability that they are all present in $F$
is $(1-\lambda n^{-1/3})^{m_{II}-1}/(2n)^{m_{II}-1}$.

\newcommand{\num}{\text{\#}}

Combining, we obtain that the expected number of bicycles $X'$
containing exactly $r$ type~I excursions, one type~II excursion,
and no type~III excursions is
\begin{align*}
\num(r,0,1)
 & \leq
\sum_{m_{II},m_1,\ldots,m_r\geq 2}
 (1/\lambda^{2r+2})n^{2r+2\over 3} (2n)^{m_{II}-2+\sum_im_i-2r}
 \notag \\ & \hspace*{1.0cm} \times
 \frac{(1-\lambda n^{-1/3})^{m_{II}-1+\sum_im_i-r}}
      {(2n)^{m_{II}-1+\sum_im_i-r}}
 \\ &=
\frac{1}{2^{r+1}\lambda^{2r+2}n^{r+1\over 3}}
 \sum_{m_{II},m_1,\ldots,m_r\geq 2}
  (1-\lambda n^{-1/3})^{m_{II}+\sum_im_i-r-1} .
\end{align*}
Note that
\begin{align*}
\lefteqn{
\sum_{m_{II},m_1,\ldots,m_r\geq 2}
 (1-\lambda n^{-1/3})^{m_{II}+\sum_im_i-r-1}
}
\\ &=
\sum_{m_{II},m_1,\ldots,m_r\geq 1}
 (1-\lambda n^{-1/3})^{m_{II}+\sum_im_i}
\\ &=
\left(\sum_{m\geq 1}(1-\lambda n^{-1/3})^m\right)^{r+1}
\\ & \leq
\frac{1}{\lambda^{r+1} n^{-{r+1\over 3}}} .
\end{align*}
Applying this to the equality above we obtain
\begin{align*}
\num(r,0,1)
 & \leq
 \frac{1}{2^{r+1}\lambda^{3r+3}} \text{, and}
\\
\sum_{r \geq 0} \num(r,0,1)
 & \leq \frac{1}{2\lambda^3(1-2\lambda^3)}
 = O(1/\lambda^3) .
\end{align*}

With similar calculations for $\num(r,\cdot,\cdot)$ this
establishes \eqref{conditioning},
and completes the proof of the case $\l \leq -1$ of
Theorem~\ref{MainScalingWindow}.
\qed

\subsection{Case $c=1+\lambda n^{-1/3}$, $\lambda\geq 1$}
\label{sparsewindow}

The proof of this part resembles the alternate proof of
Theorem~\ref{MainLowDens}.
There we showed that $m_1(t)$ a.s.\ a.e.\ followed a parabolic trajectory.
Both there and here, at time $t=\e n$,
the expectation given by the parabola is $\frac12 \e^2 n$,
and the typical deviations (the standard deviation)
from summing $\e n$ binomial r.v.s with
distributions near to $B(n,1/n)$ is about $\sqrt{\e n}$.

\begin{figure}[htbp]
\begin{center}
\includegraphics[width=4.5in,height=2.0in]{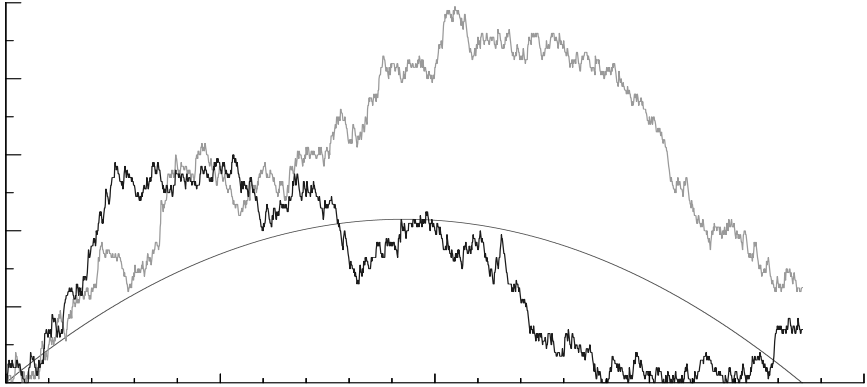}
\end{center}
\caption{Nominal parabolic trajectory of $m_1(t)$ vs $t$,
and two random samples for density $1+\l n^{-1/3}$
($\l = 2$, $n=10,000$).
With density $1+\Theta(n^{-1/3})$, the random fluctuations are
of the same order as the nominal values.}
\label{trajectories}
\end{figure}

In the previous case, with $\e=\Theta(1)$, the deviations were a.s.\
tiny compared with the expectation,
but here, with $\e = \l n^{-1/3}$,
the standard deviation of $\sqrt{\l} n^{1/3}$ is of the same order
(in terms of $n$)
as the expectation of $\frac12 \l^2 n^{1/3}$:
the trajectory is not predictable in an a.s.\ a.e.\ sense.
Figure~\ref{trajectories} shows two typical samples
(with $\l=2$ and $n=10,000$)
against the nominal parabolic trajectory.
The analysis is thus more involved.

As before, we analyze the unit-clause resolution algorithm in
which if there are any unit clauses (if $m_1(t)>0$) we choose one
at random and set its literal True,
and otherwise we choose a random literal
(from the variables not already set)
and set it True.

Our analysis proceeds in three phases.
Phase~I proceeds until time $T=2 \e n$, and we show that in this period,
there is an exponentially small chance that $m_1$ is ever
much larger than its expectation.
In Phase~II, we continue unit-clause resolution
until $m_1(t)=0$; we show that this
happens quickly.
These will give the required bounds on the integrated number of
unit clauses, and in turn unsatisfied clauses,
produced by unit-clause resolution.
In Phase~III we have a formula of density $\leq 1-\e n$, and we simply
apply the (non-algorithmic) proof of the Theorem's case $\l \leq -1$,
proved in Section~\ref{sparsewindow}.

\subsubsection{Useful facts}
We first establish a simple relation, useful for Phase~I and
essential for Phase~II.
The number of 2-clauses remaining (both of whose variables remain)
at time $\d n$ is $m_2(\d n) \distr B(n(1+\e), (1-\d)^2)$.
Thus for all times $t \leq \frac12 n(1+\e)$
(much longer than the times $\Theta(\e n)$ in which we are interested),
\begin{align}
\Pr\left(
    \max_{\d \leq \frac12}
     \left[ m_2(\d n) - n(1+\e)(1-\d)^2 \geq n^{3/5} \right]
    \right)
 \leq \exp(-\Theta(n^{1/5})) .
 \label{m2exact}
\end{align}
We prove \eqref{m2exact} using
the Chernoff bound that for a sum $X$ of independent 0-1 Bernoulli
random variables with parameters $p_1,\ldots,p_n$ and expectation
$\mu=\sum_{i=1}^np_i$,
\begin{equation}\label{chern1}
\Pr(X\ge \mu+ \D) \le \expfn{-\D^2/(2\mu+2\D/3)} .
\end{equation}
(See for example~\cite[Theorems~2.1 and~2.8]{JLR}.)
To establish \eqref{m2exact} we take $(1+\e)n$
i.i.d.\ Bernoullis with $p_i = (1-\d)^2$.
For any fixed $\d$ in \eqref{m2exact}
this immediately gives probability
$\exp(-\Theta(n^{6/5}/ n))$,
and the sum over the $\Theta(n)$ possible values of $\d$
can be subsumed into the exponential.

In the main we will therefore assume that
\begin{align}
 m_2(\d n) \leq n(1+\e)(1-\d)^2 + n^{3/5} ,
\label{m2}
\end{align}
and deal with the failure case only at the end.

We will also need two simple distributional inequalities.
First, a Bernoulli random variable is stochastically dominated
by a similar Poisson random variable,
\begin{align*}
\Be(p) \sleq \Po(-\ln(1-p)) ,
\end{align*}
as they give equal probability to 0, and the Bernoulli's remaining
probability is entirely on 1 whereas the Poisson's is on 1 and larger values.
(Here we have written $\Be(p)$ and $\Po(-\ln(1-p))$ where
we really mean random variables with those distributions;
we shall continue this practice where convenient.)
Summing $n$ independent copies of such random variables shows that
a binomial is dominated by a similar Poisson,
\begin{align*}
B(n,p) \sleq \Po(-n \ln(1-p)) .
\end{align*}
In particular, for any $a,b=\Theta(1)$,
\begin{align}
B(an, b/n) & \sleq \Po(-an \ln(1-b/n))
 = \Po(ab+O(1/n))) .
\label{Po}
\end{align}

We also recall that the exponential moments of a Poisson random variable are
\begin{align}
\E z^{\Po(d)}
 &= \exp((z-1)d) .
\label{exppois}
\end{align}

We now analyze the unit-clause algorithm in Phases I and~II.

\subsubsection{Phase~I}

During Phase~I, assuming~\eqref{m2}, at times $t=\d n$,
\begin{align*}
m_2(t) = n(1+\e)(1-\d)^2 + O(n^{3/5}) \leq n(1+1.01 \e-2\d) ,
\end{align*}
using $\e \geq n^{-1/3}$.
Meanwhile the number of unset variables is $m(t)= n(1-\d)$,
so in particular,
\begin{align}
m_2(t)/m(t) \leq 1+1.05 \e .
\label{density1}
\end{align}

With the random variables below all independent,
the unit-clause algorithm gives
\begin{align*}
\lefteqn{
m_1(t) - m_1(t-1)
} \\
 &=
 - 1 + \one(m_1(t-1)=0) - B(m_1(t-1),1/(2m(t-1)))
   \\ & \White{=} + B(m_2(t-1), 1/(m(t-1)))
 \\ & \leq
 - 1 + \one(m_1(t-1)=0) + B(m_2(t-1), 1/(m(t-1)))
 \\ & \leq
 - 1 + \one(m_1(t-1)=0) + \Po(1+1.1 \e) ,
\end{align*}
where the last inequality uses~\eqref{Po}, \eqref{density1},
and $0.1\e \gg 1/n$.

It is easy to see that, starting from $m_1(0)=m'_1(0)$,
if
$X(t) \leq Y(t) \text{ for all $t$}$, if
\begin{align*}
m_1(t) - m_1(t-1) &= \one(m_1(t-1)=0) + X(t-1) \text{ and}
\\
m'_1(t)-m'_1(t-1) &= \one(m'_1(t-1)=0) + Y(t-1) ,
\\ \intertext{then for all $t$,}
m_1(t) &\leq m'_1(t) .
\end{align*}
(An easy proof is inductive.
The $\one(\cdot)$ term may contribute to $m_1$
and not to $m'_1$ if $m_1<m'_1$, but in that case,
the inequality still holds.)
In a similar setup but with $X(t) \sleq Y(t)$,
coupling shows that $m_1(t) \sleq m'_1(t)$.

Thus $m_1(t) \sleq m'_1(t)$ where $m'_1(0)=0$ and
\begin{align*}
m'_1(t) - m'_1(t-1)
 &=
 - 1 + \one(m'_1(t-1)=0) + \Po(1+1.1\e) .
\end{align*}
Now, let $U(t)$ be a random walk with $U(0)=0$ and independent increments
\begin{align}
U(t)-U(t-1) &=
-1+\Po(1+1.1\e) ,
\label{U}
\end{align}
and let $V(t)$ count the ``record minima'' of $U$,
so $V(0)=0$ and $V(t)=V(t-1)$
except that if $U(t) < \min_{\t < t} U(\t)$, then
$V(t)=V(t-1)+1$.
Observe that
\begin{align}
m_1(t) \sleq m'_1(t) = U(t) + V(t-1) .
\label{stochm1}
\end{align}
($V(t)$ precisely takes care of the $\one(\cdot)$ terms.)

At this point, we have reduced the behavior of the
number of unit clauses $m_1(t)$ to properties of a simple
Poisson-incremented random walk.

\subsubsubsection{Renewal process $V$}

We first dispense with $V$, by showing that
\newcommand{\Vinf}{V(\infty)}
\begin{align}
\Vinf \eqdef \sup_{t \geq 0} V(t) \sleq G( 2 \e) ,
\label{Vbound}
\end{align}
where $G(p)$ indicates a geometric random variable with parameter~$p$.
Starting from any time $t_0$ at which $U(t_0)$ is a record minimum
(at which $V(t_0)=V(t_0-1)+1$),
define $U'(\t) = U(t_0+\t)-U(t_0)+1$.
Observe that $U'(0)=1$, and the first time $\t$ for which $U(\t)=0$
gives the next time $t_0+\t$ for which $V(t_0+\t)=1$.
Thus the number of ``restarts'' of the process $U'$ is $\Vinf$.

$U'$ may be viewed as a Galton-Watson branching process observed
each time an individual gives birth
(adding $\Po(\cdot)$ offspring to the population)
and itself dies (adding $-1$).
As a super-critical Galton-Watson branching process, $U'$ has
a positive probability of non-extinction, and thus the number of
restarts (following extinctions) is geometrically distributed.

Quantitatively, the extinction probability of a Galton-Watson
process with $X$ offspring
(the probability the process never hits~0)
is well known to be
the unique root $p \in [0,1)$ of
\begin{align}
p = \E (p^X) .
\label{gw}
\end{align}
(See for example~\cite[pp.~247--248]{Durrett}.)
Also, for any $p$ such that $p > \E (p^X)$,
the probability of non-extinction exceeds $1-p$.
In this case, recalling \eqref{U} and \eqref{exppois},
we seek $p$ such that
\begin{align*}
 p &> \E(p^X) = \exp((p-1)(1+1.1\e))
\end{align*}
or equivalently, with $q=1-p$,
\begin{align*}
 \ln(1-q) &> -q (1+1.1\e) .
\end{align*}
Taking a Taylor expansion around $q=0$ and cancelling like terms,
it suffices to ensure that
$\frac12 q + \frac13 q^2 + \cdots < 1.1 \e$,
and $q= 2\e$ suffices
(for all $\e < 0.37$, let alone the $\e=\Theta(n^{-1/3})$ of interest).

Thus $U'$ has non-extinction probability at least $2 \e$,
verifying~\eqref{Vbound}.

\subsubsubsection{Random walk $U$}
\label{sec:U}

We now analyze the random walk~$U$ to show that
for any
$0< \e \leq 0.02$ and $0 < \a \leq 0.06$
(our principal realm of interest will be $\e,\a = \Theta(n^{-1/3})$),
for any time $t$,
\begin{align}
\Pr \left( \max_{0 \leq \t \leq t} U(\t)
          \geq \E U(t) + \a t \right)
 &\leq
 \exp(-t \a^2 / 2.1) .
\label{stochu}
\end{align}
Observe that $U(t)$ is a submartingale,
and for any $\b>0$
(by convexity of $\exp(\b u)$),
$\exp(\b U(t))$
is a non-negative submartingale.
It follows from Doob's submartingale inequality (see~\cite{Durrett}) that
\begin{align}
\lefteqn{
\Pr \left( \max_{0 \leq \t \leq t} U(\t)
          \geq \E U(t) + \a t \right)
} \notag \\
 &=
\Pr \left( \max_{0 \leq \t \leq t}
          \exp\left( \b U(\t) \right)
          \geq \exp\left( \b (\E U(t) + \a t \right) ) \right)
\notag
 \\ & \leq
\frac{\E \left( \exp(\b U(t)) \right)}
     {\exp\left( \b (\E U(t) + \a t) \right)} .
\label{doob}
\end{align}
Trivially,
\begin{align}
\E U(t) = -t + (1+1.1\e) t = 1.1 \e t ,
\label{eu}
\end{align}
and, by~\eqref{exppois},
\begin{align*}
\E \left( \exp(\b U(t)) \right)
 &=
\exp(-\b t + \b \Po((1+1.1\e)t))
\\ &=
\exp(-\b t) \exp((e^\b-1)(1+1.1\e)t) ,
\end{align*}
so \eqref{doob} is
\begin{align}
\exp(-t [\b - (1+1.1\e)(e^\b-1) + \b(1.1\e+\a)]) .
\label{beta}
\end{align}
We are free to choose $\b>0$ as we like, so to minimize \eqref{beta}
we maximize the innermost quantity.
Setting its derivative equal to~0 yields
$1 - (1+1.1\e) e^\b + 1.1\e + \a = 0$ or
$\b = \ln(1+\a/(1+1.1\e))$,
but we will simply take $\b=\a$.
Then (eschewing asymptotes in favor of absolute bounds),
for $\e < 0.02$ and $\a < 0.06$
(let alone the regime $\e,\a=\Theta(n^{-1/3})$ of interest),
\eqref{beta} is
\begin{align*}
 \leq \expfn{-t \a^2 / 2.1} ,
\end{align*}
proving~\eqref{stochu}.

\subsubsubsection{Parameter substitution}
\label{sec:Uparms}

In the case of interest, $1\leq \l \ll n^{1/3}$, $\e=\l n^{-1/3} \ll 1$,
and $t = 2 \e n = 2 \l n^{2/3}$.
Here, \
\begin{align*}
\E m_1(t)
 \leq \E U(t) + \E \Vinf
 \leq 1.1 \e t + \frac1{2\e}
 \leq 2.2 \l^2 n^{1/3} + n^{1/3}
 \leq 3.2 \l^2 n^{1/3} .
\end{align*}
Substituting $\a= \a' / \sqrt{t}$ into \eqref{stochu}
(which is then valid up to $\a' = n^{1/4} = o(n^{1/3})$)
\begin{align}
\Pr(\max_{\t \leq t} U(t) \geq \E U(t) + \a' \sqrt{2 \l} n^{1/3})
 \leq \exp(-\a'^2 / 2.1) ,
\label{utail}
\end{align}
so the tails of $U(t)$ fall off exponentially with a ``half-life''
smaller than the bound on the mean
(as $\l \geq 1$ implies $\sqrt{2 \l} < 2.2 \l^2$).
$\Vinf$ has an expectation which is at most comparable,
and (as a geometric random variable) again falls off
exponentially with half-life comparable to its mean.

It follows that
\begin{align*}
\E \sum_{\t=1}^{2 \e n} m_1(\t)
 & \leq (2 \e n) \E(U(2\e n) + \Vinf)
 \leq 6.4 \l^3 n,
\end{align*}
and,
for $\a' \leq n^{1/4}$,
\begin{align}
\Pr \left( \max_{\t \leq 2 \e n} m_1(\t)
 \geq \a' 3.2 \l^2 n^{1/3} \right)
 &= \exp(-\Omega(\a'))
 \label{m1final}
 \text{ and}
 \\
\Pr \left( \sum_{\t=1}^{2 \e n} m_1(\t)
 \geq \a' 6.4 \l^3 n \right)
 &= \exp(-\Omega(\a')) .
 \label{m1area}
\end{align}

The probability of a deviation with $\a' > n^{1/4}$
is $\exp(-\Omega(n^{1/2}))$, and will be dealt with
as a ``failure probability'' at the end.

\subsubsection{Phase II}

The analysis of this phase largely parallels the previous one.

Assuming~\eqref{m2}, at times $t=\d n$,
$m_2(t)/m(t)$ is roughly $(1+\e)(1-\d)$, and in particular,
since in Phase~II by definition $\d \geq 2 \e$,
\begin{align}
m_2(t)/m(t) \leq 1- 0.95 \e .
\label{density2}
\end{align}

Since Phase~II ends as soon as $m_1(t) = 0$,
there is no $+\one(\cdot)$ term to worry about,
so assuming~\eqref{density2},
\begin{align*}
\lefteqn{
m_1(t) - m_1(t-1)
} \\
 &=
 - 1 - B(m_1(t-1),1/(2m(t-1)))
   + B(m_2(t-1), 1/(m(t-1)))
 \\ & \leq
 - 1 + \Po(1 - 0.9 \e) .
\end{align*}
By the same argument as for Phase~I, then,
\begin{align*}
m_1(2\e n + t) \sleq m_1(2 \e n) + W(t)
\end{align*}
where $W(t)$ is a random walk with $W(0)=0$ and
independent increments $-1+\Po(1-0.9\e)$.

Note that $W(\e n) = W(\l n^{2/3})$ has mean and standard deviation
both $\Theta(n^{1/3})$, so for multiples of this time,
$W$ is exponentially sure to achieve at least half its
(negative) expectation;
we now quantify this.
At time $\a \e n$,
\begin{align*}
\lefteqn{
\Pr(W(\a \e n) > -\tfrac12 \a \e^2 n)
} \\
 &=
 \Pr\left( \Po(\a \e n (1-.9\e)) > \E(\Po(\cdot)) + 0.4 \a \e^2 n \right)
 \notag
 \\ & \leq
 \expfn{-\frac{(0.4 \a \e^2 n)^2}{\a \e n(1-0.9\e) + 0.4 \a \e^2 n}}
 \notag
\end{align*}
since the Chernoff bound~\eqref{chern1} applies
as well to the Poisson.
Substituting $\e=\l n^{-1/3}$, the
denominator's first term, of order $\Theta(\a n^{2/3})$,
dominates the second, of order $\Theta(\a n^{1/3})$, giving
\begin{align}
\Pr(W(\a \e n) > -\tfrac12 \a \e^2 n)
 & \leq
 \exp(-0.4^2 \a \l^3) .
 \label{2end}
\end{align}

Then, conditionally on Phase~I ending at $m_1(2\e n) = \a_1 \l^2 n^{1/3}$
(see~\eqref{m1final}),
for any $\a > 2 \a_1$,
\eqref{2end} implies that
Phase~II ends by time $2\e n+\a \e n$,
with probability exponential in $\a$.

Furthermore, over Phase~II, $m_1(t)$ is unlikely ever to increase
much over its initial value.
An argument along the lines used in the context of
equation~\eqref{gw} could be constructed to show that
$\max_{t \geq 0} W(t)$ is exponentially sure to be quite small,
but as there are some technical complications,
we take a simple, wasteful approach.
Observe that
\begin{align*}
W(t) \sleq X(t)
\end{align*}
where $X(0)=W(0)=0$ and $X(t)$ has independent increments $-1+\Po(1+\e)$.
This wild over-estimation is useful because $X$ (unlike $W$) is a
submartingale, to which we apply Doob's inequality.
Just as in sections \ref{sec:U} and~\ref{sec:Uparms},
over a time interval $\a \e n$,
$X$ is exponentially unlikely ever to exceed
a multiple of its final expectation,
$\E X(\a \e n) = \a \e^2 n = \a \l^2 n^{1/3}$,
and so $W$ and in turn $m_1$ are
at least as unlikely to rise more than this
amount above their initial values.

So, conditionally on Phase~II starting at $m_1 = \a_1 \l^2 n^{1/3}$,
Phase~II finishes within additional time $\a \l^2 n^{1/3}$
with exponentially high probability for all $\a \geq 2 \a_1$,
and within that additional time,
$m_1$ is exponentially unlikely (in $\b$) to exceed
$(\a_1 + \b \a) \l^2 n^{1/3}$.
It follows that if Phase~II ends at time $2\e n+t$,
\begin{align}
\Pr\left(
\sum_{\t=1}^{t} m_1(2\e n+\t)
 \geq
 \const \b' a_1^2 \l^3 n
 \right)
 & \leq \exp(-\b') .
 \label{phaseIIarea}
\end{align}

\subsubsection{Phases I, II and III}
We have argued that over Phases~I and~II
the number of unit clauses $m_1(t)$ is exponentially
unlikely ever to exceed a multiple of $\e^2 n = \l^2 n^{1/3}$,
and that Phase~II is exponentially unlikely to end after
a multiple of time $\e n = \l n^{2/3}$,
to prove, in~\eqref{m1area} and~\eqref{phaseIIarea},
that the summed number of unit clauses
$M_1 = \sum_{\t} m_1(\t)$
(summed over times $\t$ from 0 to the end of phase~II), is exponentially
unlikely
to exceed a multiple of $\l^3 n$:
\begin{align*}
\Pr(M_1 \geq \const \b \l^3 n) \leq \exp(-\b) .
\end{align*}

By definition of the unit-clause algorithm, at each stage the
literals forming the unit clauses are
drawn independently at random with replacement
from among the literals not yet set,
and so the number of unit clauses dissatisfied at each step $t$
is
\begin{align}
B(m_1(t), 1/(2 (n-t))
\label{m1unsat}
\end{align}
(where $m_1(t)$ is itself a random variable).
With probability $1-\exp(-\Theta(n^{1/4}))$ these phases end
long before time $t=n/3$,
so \eqref{m1unsat} is $\sleq \Po(0.8 m_1(t)/n)$,
and by independence of the random variables in \eqref{m1unsat}
(each conditioned on $m_1(t)$)
for different times~$t$,
the total number of unit clauses dissatisfied in phases~I and~II
is dominated by $\Po(0.8 M_1/n)$.

Since $\E M_1 = O(\l^3 n)$,
the Poisson's expectation is $O(\l^3)$, and
the number of $X$ of unit clauses unsatisfied over these phases
also has $\E X = O(\l^3)$;
this confirms (for Phases~I and~II) \emph{one} assertion of
Theorem~\ref{MainScalingWindow}).
Fixing $\b=1$, there is at least constant probability that
$M_1 \leq \const \l^3 n$ and so
the probability that \emph{no} unit clause is dissatisfied is
$\Pr(X=0) \geq \exp(-O(\l^3))$, a \emph{second} assertion of the theorem.
Since both $M_1$ and $\Po(M_1/n)$ have exponential tails, so does $X$
---
$\Pr(X \geq \b \const \l^3) \leq \exp(-\b)$ ---
a \emph{third} assertion of Theorem~\ref{MainScalingWindow}.
We now argue that Phase~III leaves all these properties intact.

By construction, at the conclusion of Phases~I and~II the
remaining formula is uniformly random, still on $n(1-o(1))$
variables, but now with density $\leq 1-\e \leq 1-n^{-1/3}$.
For Phase~III we simply argue that,
by the previously proved case $\l \leq -1$ of this Theorem,
such a formula can be satisfied but for $\leq \const \b$ clauses,
with probability $\geq 1-\exp(-\b)$.
This concludes the proof of the case $\l>1$
of Theorem~\ref{MainScalingWindow}.
\qed

\subsubsection{Remarks}

A corresponding \emph{lower bound}, on the number of clauses that
must be violated, cannot be found by the same techniques, since
there is no guarantee that the unit-clause algorithm is doing the
best possible.
One alternative is to analyze the pure-literal rule, which is
guaranteed to make no ``mistakes'' as long as it runs,
then use other methods to analyze the remaining ``core'' formula;
we understand that this analysis has been done successfully
(and independently) by Kim~\cite{KimPoissonCloning02}.
Another approach might be to extend the ``bicycles'' analysis of
Theorem~\ref{subthreshold}
(or the $\l \leq -1$ case of Theorem~\ref{MainScalingWindow})
to the case~$\e>0$ (particularly, $\e=\l n^{-1/3}$ and $\l>1$),
but this seems not to be easy.

We remark that,
no matter the particular approach pursued, verifying that
the number of clauses that must be dissatisfied is $\Omega(\l^3)$
seems to lead back, in intuition and in proof techniques,
to the fact that in a $G_{n,p}$ random graph with average degree
$np = 1+\l n^{-1/3}$,
there is likely to be a giant component whose ``kernel''
is a random cubic graph on $\Theta(\l^3)$ vertices~\cite[p.~123]{JLR}.

\section{Random MAX $\mathbf{k}$-SAT and MAX CSP}
\label{sec:ksat}

In this section we present some general facts and conjectures
about \maxksat and \maxcsp, and
generalize the 2-\sat high-density results.

\subsection{Concentration and limits}
\label{sec:concentration}

It is known that random $k$-\sat has a sharp threshold:
that is,
there exists a threshold function $c(n)$
such that for any $\e>0$,
as $n \ra \infty$,
a random formula on $n$ variables
with $(c(n)-\e)n$ clauses is a.s.\ satisfiable,
while one
with $(c(n)+\e)n$ clauses is a.s.\ unsatisfiable~\cite{Frie}.
To prove an analogous result for random \maxksat is much easier;
this was first done by~\cite{BFU93}.
We will
employ a ``bounded difference'' inequality; specifically, a
generalization of Azuma's inequality in a form due to
McDiarmid~\cite{CMD89} (see also Bollob\'as~\cite{Bo88}).

\begin{theorem}[Azuma]
\label{thm:azuma}
Let $X_1, \ldots, X_n$ be independent random variables,
with $X_k$ taking values in a set $A_k$ for each $k$.
Suppose that the (measurable) function $f: \prod A_k \ra \Re$ satisfies
$|f(x) - f(x')| \leq c_k$
whenever the vectors $x$ and $x'$ differ only in the $k$'th coordinate.
Let $Y$ be the random variable $f(X_1, \ldots, X_n)$.
Then for any $\l>0$,
$ \pr[\,|Y-\E Y| \geq \l] \leq 2 \exp(-2 \l^2/ \sum c_k^2) $.
\end{theorem}

Let $F_k(n,m)$ be a random $k$-\sat formula on $n$ variables
with $m$ clauses, and
let $f_k(n,m) = \E(\max F_k)$;
we may omit the subscripts $k$.

\begin{theorem}\textnormal{(\cite{BFU93})}
\label{concentration}
For all $k$, $n$, $c$, and $\lambda$,
$\pr(| \max F_k(n,cn) - f_k(n,cn) | > \lambda)
 < 2 \exp(-2 \lambda^2 / (cn)) $.
\end{theorem}

\begin{proof}
Let $X_i$ represent the $i$th clause in~$F$.
Replacing $X_i$ with an arbitrary clause cannot
change $\max F$ by more than~1.
The result follows from Azuma's inequality.
\end{proof}

The theorem's statement that for any $c$ and large $n$,
$F(n,cn)/(cn)$ has some almost-sure almost-exact value,
is reminiscent of Friedgut's theorem (Theorem~\ref{thm:friedgut}) that
(loosely interpreted) says that
for large $n$ and any $c$ away from the threshold,
$\Pr(F(n,cn) \text{ is satisfiable})$
is almost exactly either~0 or~1.
In our case, the target value $f(n,cn)/(cn)$
is unknown and it is unknown whether it has a limit in~$n$,
and in Friedgut's case, again,
it is unknown for which $c$s
the probability is near 0 and for which it is near~1,
and whether the threshold value of $c$
(and the distribution function) has a limit in~$n$.
To conjecture that $f(n,cn)/(cn)$ tends to a limit in $n$
is in this sense
analogous to the ``satisfiability threshold conjecture''.

\begin{conjecture}\textnormal{(\maxsat limiting function conjecture)}
\label{conjecture:limit}
For every $k$,
for every constant $c>0$,
as $n \ra \infty$,
$f_k(n,cn)/n$ converges to a limit.
\end{conjecture}
The conjecture
may equally well be extended to arbitrary~\csp{s},
yet is open even for \maxtsat.

If $f_k(n,cn)/(cn)$ were monotone in $n$, the conjecture's
truth would follow.
Of course we do not know this, but can prove monotonicity
in~$c$: that as the number of clauses increases,
the expected fraction of clauses that can be satisfied can only decrease.

\begin{remark}
For any $k$ and $n$, $f_k(n,m) / m$ is a non-increasing function of~$m$.
\end{remark}

\begin{proof}
In a uniform random instance of $F_k(n,m)$, let the maximum number
of satisfiable clauses be $J$, so that $\E(J)=f(n,m)$.
By deleting single clauses,
we obtain $m$ uniform random instances $F$ of $F(n,m-1)$.
Of these, $m-J$ each have $\max F=J$, while the
remaining $J$ each have $\max F\in \{J-1,J\}$.
The average of these $m$ values is at least
$\frac{(m-J)(J)+(J)(J-1)}{m}=\frac{J(m-1)}{m}$.
Taking expectations, we find
$\frac{f(n,m-1)}{m-1} \geq
 \frac{1}{m-1} \times \E(\frac{J(m-1)}{m}) =
 \E(\frac{J}{m}) = \frac{f(n,m)}{m}$,
as desired.
\end{proof}

Finally, we expect a connection between the \maxsat limiting
function conjecture (Conjecture~\ref{conjecture:limit}) above
and the usual satisfiability threshold conjecture
(Conjecture~\ref{conj:thresh}).
We formalize this in the following conjecture.

\begin{conjecture}
\label{conj2}
For any $c<0$, $\lim_{n \ra \infty} f(n,cn)/(cn) = 1$
if and only if $\lim_{n \ra \infty} \Pr(F(n,cn) \text{is satisfiable}) = 1$.
\end{conjecture}

One aspect of this is easily resolved.
If $\limsup f(n,cn)/(cn) < 1$, say $1-\d$,
then on average $c \d n$ clauses per formula go unsatisfied,
at least a $\d$ fraction of all formulas must be unsatisfiable,
and so $\limsup \Pr(F(n,cn) \text{ is satisfiable}) < 1$.
But nothing more seems obvious.

\subsection{High-density MAX $\mathbf{k}$-SAT and MAX CSP}
In this section we extend Theorem~\ref{MainHighDens}.

\begin{theorem}
\label{claim2}
\label{ksat1}
For all $k$, for all $c$ sufficiently large,
$ (\frac{2^k-1}{2^k}c +
 \frac{2}{k+1} \sqrt{\frac{c k}{\pi 2^k}} - O(1) )n
\leqa f_k(n,cn) \leqa
 (\frac{2^k-1}{2^k}c+ \sqrt{c} \sqrt{\frac{(2^k-1) \ln 2}{2^{2k-1}}}) n $.
\end{theorem}

\noindent Note that the leading terms are equal, and the
second-order terms equal to within $\const \cdot \sqrt{k}$.

\begin{proof}
\textbf{Upper bound.}
The proof is very similar to that of Theorem~\ref{MainHighDens}.
Using the first-moment method, we have:
\begin{align*}
P &= \pr(\exists \text{ satisfiable } F')
 \\ & \leq
 2^n \sum_{l=0}^{r c n} \binom{cn}{l}
  (\frac{2^k- 1}{2^k})^{cn-l} (\frac{1}{2^k})^l.
\end{align*}
For $r<\frac{1}{2^k}$ the sum is dominated by the last term, and
so we fix $l=r c n$.  Using (\ref{binom}), taking logarithms, and
finally substituting $r= \frac{1}{2^k}-\e$, we have
\begin{align*}
\frac{1}{cn}\ln P \asymp
\frac{\ln(2)}{c}-(\frac{2^{2k-1}}{2^k-1})\e^2+ O(\e^3).
\end{align*}
Thus for $r<1/(2^k)-\sqrt{\frac{(2^k-1)\ln 2}{c2^{2k-1}}}$, $P\ra 0$ as $n \ra \infty$.

\textbf{Lower bound.}
Set the variables sequentially.
Set variables $X_1,X_2,\ldots,X_{\ell-1}$ randomly,
and then
\newcommand{\notell}{i}
for each $\ell \leq \notell \leq n$,
enumerate those clauses involving only $X_\notell$ and some subset of
$\{X_1,X_2,\ldots,X_{\ell-1}\}$ (that is, unit clauses).
The expected number of such clauses is about
$$c n k ( \frac{1}{n})(\frac{\ell-1}{n})^{k-1} =
 c k (\frac{\ell-1}{n})^{k-1},$$
and if we count only those left unsatisfied by their previous $k-1$
variables, the expected number becomes
$$h_\ell=\frac{c k}{2^{k-1}}(\frac{\ell-1}{n})^{k-1}.$$
(Here we incur a minor error by sampling with replacement instead of
without; $(\frac{\ell-1}{n})^{k-1}$ should really be
$\prod_{1 \leq h \leq k-1}(\frac{\ell-h}{n-h})$.)
More precisely, the number of such clauses enjoys a Poisson distribution
with mean $h\ell$.
Set the value of $X_\notell$ to maximize the number of such clauses
satisfied; as before, this number is about
$\frac{1}{2}h_\ell + \frac{1}{2} \sqrt {\frac{2}{\pi} h_\ell} + O(1)$.
The advantage over purely random guessing is
$$\sqrt{\frac{1}{2\pi} h_\ell} + O(1)=
\sqrt {\frac{c k}{2 \pi 2^{k-1}}(1-\frac{\ell-1}{n})^{k-1}}+O(1).$$
Sum over
$\notell = \ell, \ldots, n$
to obtain an advantage of
$$\sqrt{\frac{c k}{\pi 2^k}} \frac{2n}{k+1}+O(n).$$

\end{proof}

Still more generally, we may consider a \csp
(constraint satisfaction problem).
Let $g$ be a $k$-ary ``constraint'' function,
$g: \{0,1\}^k \ra \{0,1\}$.
A random formula $F_g(n,m)$ over $g$ is defined by
$m$ clauses, each chosen uniformly at random (with replacement)
from the $2^k n (n-1) \cdots (n-k+1)$ possible clauses
defined by an ordered $k$-tuple of distinct variables
each appearing positively or negated.
(Formally, a clause consists of a $k$-tuple
$(i_1,\ldots,i_k)$ of distinct values in $[n]$, specifying the variables,
and a binary $k$-vector $(\sigma_1,\ldots,\sigma_k)$, specifying their signs.)
A clause with variables (signed variables) $X_1,\ldots,X_k$
is satisfied if $g(X_1,\ldots,X_k)=1$.
(Formally, an assignment $x_1,\ldots,x_n$
of the full set of variables $X_1,\ldots,X_n$
satisfies a clause as above
if $g(x_{i_1} \oplus \sigma_1, \ldots, x_{i_n} \oplus \sigma_n) = 1$,
where ``$\oplus$'' denotes \xor, or addition modulo~2.)
As ever, such a formula $F$ is satisfiable if there exists
an assignment of the variables satisfying all the clauses;
and $\max F$ is the maximum, over all assignments,
of the number of clauses satisfied.

Generally a \csp may be based on a finite family of
constraint functions, of ``arities'' bounded by~$k$,
but for notational convenience we limit ourselves to a single function.

Let a $k$-ary clause function $g$ be given, with $\E(g(X)) = p$
over random inputs.
Define $P=\min\{p,1-p\}$ and $Q=1-P$.
Let $F_g(n,m)$ be a random formula over $g$ on $n$ variables,
with $m$ clauses,
and let $f_g(n,m) = \E(\max F)$.

\begin{theorem}
\label{csp}
Given an arity $k$ and a constraint function $g$,
for all $c$ sufficiently large,
$
(p \, c + \sqrt{P Q^2 c/k} \, ) n
\leqa f_g(n,m) \leqa
(p \, c + \sqrt{2PQ \ln(2) c} \, )n
$.
\end{theorem}

The proof follows that of Theorem~\ref{ksat1}, and is omitted.

\section{Online random MAX 2-SAT}
\label{sec:online} In this section, we discuss online versions of
the \maxtsat problem.
\cite{BF01,BFW02} consider an online version of \mgfss,
in which random edges $e_i$ are given one by one,
and we must accept or reject $e_i$
based on the previous edges $e_1,\ldots,e_{i-1}$,
with the goal of accepting as many edges as possible
without creating a giant component.

There are two natural online interpretations of
random \maxtsat.
In both,
we are told in advance the total number of variables $n$
and clauses~$m$;
also, in both,
clauses $c_i$ are presented one by one,
and we must choose ``on line'' whether to accept or reject $c_i$
based on the previously seen clauses $c_1, \ldots, c_{i-1}$.
When we accept a clause we are guaranteeing to satisfy it;
when we reject a clause we are free to satisfy or dissatisfy it.
Our goal is to maximize the number of clauses accepted.

In our first interpretation of online \maxtsat, \problema,
when we accept a clause,
we are also required to satisfy it immediately,
by setting at least one of its literals True;
once a variable is set, it may never be changed.
The second interpretation, \problemb, is more generous:
the variables' assignments may be decided
after the last clause is presented.
Let $\foa(n,m)$ be the expected number of clauses accepted by
an optimal algorithm for \problema,
and $\fob(n,m)$ that for \problemb.
Clearly,
$\frac34 m \leq \foa(n,m) \leq \fob(n,m) \leq f(n,m)$.
Here we present a ``lazy'' algorithm applicable to
both $\foa(n,cn)$ and $\fob(n,cn)$.
\onlinea
begins with no variables ``set''.
On presentation of a clause, \onlinea rejects it only if it must,
and otherwise does the least it can to accept it.
Specifically,
on presentation of clause $c_i$, which without loss of generality
we may consider to be $(X \vee Y)$, it takes the following action.
If $X=\true$ or $Y=\true$, accept $c_i$.
If $X=\false$ and $Y=\false$, reject $c_i$.
If $X=\false$ and $Y$ is unset (or vice-versa),
set $Y=\true$ (resp.\ $X=\true$) and accept $c_i$.
If $X$ and $Y$ are both unset,
arbitrarily choose one, set it \true, and accept $c_i$.

\begin{theorem}
\label{online}
For any fixed $c$,
\onlinea is the unique (up to its arbitrary choice)
optimal algorithm for \problema, and
$\foa(n,cn) \asymp (\frac34 c + (1-e^{-c})/4 + (1-e^{-c})^2/8) n
 \geq (\frac34 c + \frac38) n$.
\end{theorem}

\noindent
We note that for $c=1$,
$\foa(n,n) \approx 0.957997 n$,
and for $c$ asymptotically large,
$\foa(n,cn) \asymp (\frac34 c + \frac38) n$.

\pf{Proof of optimality.}
On appearance of a clause $c_i$, it is clearly best not to set any
variable not appearing in $c_i$, for this merely imposes extra constraints.
Similarly, if $c_i$ is already satisfied by one of its literals,
then it is best to accept it and to set no additional variables.

The only interesting cases, then, are if $c_i$ is not already satisfied,
but one or both of its variables are unset.
Again, if both variables are unset, it is best to set at most one of them,
and it doesn't matter which one:
the ``future'' performance of an optimal algorithm is solely
a (random) function of the number of unset variables
and the number of clauses remaining,
and these parameters of the future, as well as the number of clauses
accepted in the past, are the same whether $c_i$'s first or second
literal is set.

It only remains to show that if $c_i$ is not satisfied by a variable
already set, and at least one of its variables is not yet set,
then an optimal algorithm must set one of its literals to \true.
Consider a putatively optimal algorithm \Opt which does not do this,
so for a literal $X$ in $c_i$, either \Opt sets $X$ to \false,
or it leaves $X$ unset.

In the case when \Opt sets $X$ to \false,
let a competing algorithm $\Opt'$ set $X$ to \true,
then simulate \Opt but reversing the roles of $X$ and $\Xbar$
in future clauses.
``Couple'' the distribution of future random clauses seen by
\Opt and $\Opt'$, also by reversing the roles of $X$ and~$\Xbar$.
With this coupling, $\Opt'$ accepts
exactly the same number of clauses as $\Opt$ in the future,
but has accepted one additional clause so far ($c_i$);
this contradicts the supposed optimality of~\Opt.

The slightly less obvious case is when \Opt leaves $X$ unset.
Again we introduce a competing algorithm $\Opt'$,
which sets $X$ to \true, then simulates \Opt until such time
as \Opt sets $X$.
For inputs where \Opt never sets $X$, $\Opt'$ accepts every
clause that $\Opt$ accepts, as well as the clause $c_i$,
and perhaps additional clauses in which $X$ appears;
$\Opt'$ is strictly better on these inputs.
For inputs where \Opt eventually sets $X$ to \true, $\Opt'$ goes
on simulating $\Opt$, again peforming exactly as well on future
clauses, and strictly better on past ones.
For inputs where at time $j>i$,
\Opt sets $X$ to \false,
$\Opt'$ may simulate $\Opt$ but (as before) with the roles
of $X$ and $\Xbar$ reversed.
With the previous coupling, on these inputs,
$\Opt'$ accepts exactly as many future clauses as $\Opt$,
and at least as many in the past
($\Opt'$ has accepted $c_i$ and perhaps other clauses rejected by $\Opt$,
while \Opt has accepted $c_j$ and no other clause rejected by $\Opt'$).
So in all three cases, the expected number of clauses accepted by
$\Opt'$ is at least as many as for \Opt,
and in the first two cases, which occur with nonzero probability
(for example, if no future clause contains $X$),
strictly more;
this contradicts the supposed optimality of \Opt.

\pf{Proof of performance.}

Note that clauses causing a variable to be set by \onlinea
are always satisfied, and
those not causing a variable to be set are satisfied with
probability $3/4$ (if both variables are set) or $1$ (if one is
set satisfyingly).

If $k$ variables are yet to be set,
the probability that a clause has neither variable set is
$(k/n)^2$, the probability it has one variable set
non-satisfyingly and the other not set is $2 \cdot \frac12 \cdot
((n-k)/n)(k/n)$, so a random clause falls into one of these cases
w.p.~$k/n$. The expected time to set another variable when $k$ are
unset is thus $n/k$.
In this period, clauses have (unconditioned) probabilities
$(n-k)^2/n^2$ that both variables are set, and $k(n-k)/n^2$ that
one is set satisfyingly and the other unset; conditional upon one
or other of these being the case (a variable is not set for this
clause), the probabilities are $(n-k)/n$ for the first case and
$k/n$ for the second, and the clause is satisfied with
probabilities $3/4$ and $1$ in these cases, for average gain
$\frac{1}{4}k/n$ over the naive $3/4$. The total
gain in the number of clauses satisfied in the expected $n/k-1$
steps before the setting, and the $n/k$'th step with the setting,
is $(\frac{n}{k}-1)(\frac14 k/n) + 1/4
 = 1/2 - \frac14 k/n$.
The process goes through $k=n, n-1, \ldots, n-I^\star$, until
the sum of the waiting times  exceeds the number of clauses $cn$.
Where $H(i)$ denotes the $i$'th harmonic number,
for a given~$I$,
the \emph{expected} sum of the
waiting times is $\sum_{i=0}^I n/(n-i)
 = n (H(n)-H(n-I-1))
 \approx n (\ln(n/(n-I)))$.
Solving for this equal to $cn$ gives $n/(n-\hat{I}) = \exp(c)$,
or $\hat{I} = n (1-\exp(-c))$.

What is the variance in the total waiting time~$W$, for $I=\hat{I}$,
and where we will allow the total to exceed~$cn$?
Each individual waiting time is geometrically distributed
with a mean in the range
$\frac{n}{n}= 1$ to $\frac{n}{n (\exp(-c))}= \exp(c)$,
all of which are $O(1)$,
so $W$ has standard deviation~$O(\sqrt{n})$.
The amount by which we may have overshot (or fallen short of)
the target value~$cn$ is $W-cn$;
since each round takes time at least~1,
to reach precisely $cn$ it suffices to back off (or add)
at most $W-cn$ rounds.
That is, $|I^\star -\hat{I}| \leq |W-cn|$,
which with probability exponentially close to~1 is $O(n^{2/3})$.
The expected total number of clauses satisfied over the naive
$3/4$ fraction is then
$\E \left( \sum_{i=0}^{I^\star} (1/2 - \frac14 (n-i)/n) \right)
 \asymp \hat{I}/4 + \hat{I}^2 / (8n)$.
That is, the expected number of clauses satisfied is
 $\asymp (\frac34 c + (1-e^{-c})/4 + (1-e^{-c})^2/8) n$.

\qed

Note that \onlinea does not, in fact, need to know the
number of clauses in advance!

A variant of \problema is that if we accept a
clause we must set \emph{both} its variables.
In this case, similar arguments show that
an optimal algorithm simply sets each new literal \true.

We know essentially nothing about \problemb.
To obtain improved bounds, or, ideally,
to identify a provably optimal algorithm,
are interesting open problems.

\section{Random MAX CUT}
\label{sec:cut}

\subsection{Motivation}
One source of motivation for our work was, as mentioned in the
introduction,
that although \emph{random} constraint satisfaction problems (\csp{s})
and \emph{max} \csp{s} are well studied, random \maxsc \csp{s}
seem not to have been.
However, we had a second, particular source of motivation,
in recent work on ``avoiding a giant component'' in a random graph.

Think of \maxsat as the problem of, given a formula,
to select as many clauses as possible so that the
subformula of selected edges is satisfiable.
An analogous problem is, given a graph,
to select as many edges as possible so that the
subgraph of selected edges has no giant component
(suitably defined).

The latter problem was posed in a slightly different form
by Achlioptas, who asked how many random edge \emph{pairs} could be
given, such that by selecting one edge from each pair,
a giant component could be avoided.
Bohman and Frieze showed in \cite{BF01} that a giant component
can be avoided with $0.55 n$ edge pairs
(where a random selection of one edge from each pair would
almost surely generate a giant component).
Bohman, Frieze, and Wormald \cite{BFW02} considered the problem
without Achlioptas's original ``pairing'' aspect:
how many edges may a random graph have,
so that some subgraph with $1/2$ the edges has no giant component.
They show that this is true up to about $1.958 n$ edges but not
beyond (where the precise threshold satisfies a transcendental equation).
Without the pairing aspect, there is no longer anything special
about $1/2$, though, and \cite{BFW02} is easily extended to
answer the question: for a random graph $G(n,cn)$,
how many edges $f(n,cn)$ may be retained while avoiding a giant component.
This is precisely the same sort of question we considered for \sat,
and was in our minds when we began this work.

It is tempting to imagine a particular connection between the
two questions, because of a well known
connection between the unsatisfiability of a random 2-\sat formula
and the existence of a giant component in a random graph,
most easily explained in terms of branching processes.
For a 2-\sat formula~$F$, consider a branching process on literals,
where a literal $X$ has offspring including $Y$ if $F$ includes
a clause $\set{\Xbar,Y}$ (and if $Y$ was not the parent of $X$).
(The process models the fact that if $X$ is set true,
$Y$ must also be set true to satisfy~$F$).
Although additional work is needed to prove it, a random 2-\sat
formula is satisfiable with high probability if this branching process
is subcritical (if each $X$ has an expected number of offspring~$<1$)
and unsatisfiable \whp if it is supercritical.
For a random graph $G$, consider a branching process on vertices,
where a vertex $v$ has offspring including $w$ if
$G$ has an edge $\set{v,w}$ (and if $w$ was not the parent of $v$).
Here, \whp $G$ has no giant component if the process is subcritical,
and \whp has one if it is supercritical.
These intuitively explain the phase-transition thresholds
of $cn$ clauses, $c=1$, for a random 2-\sat formula,
and edge density $c/n$, $c=1$, for a random graph.

Despite this connection between unsatisfiability of a random formula,
and a giant component in a random graph,
the size of a largest giant-free subgraph of a random graph
behaves very differently from the size of a largest satisfiable
subformula of a random formula.
Specifically, for large clause density $c$, there is a satisfiable subformula
preserving an expected constant fraction ($3/4$ths) of the clauses,
while for a random graph with $cn$ edges,
the largest giant-free subgraph has only about $n$ edges, a $1/c$ fraction.
This can be read off from Theorem~\ref{bfwplus},
or argued more simply:
if $G$ had a giant-free subgraph $H$ with linearly more than $n$ edges,
$H$ (and thus $G$) would have to have a linear-size dense component,
but a random sparse graph has no linear-size dense component.

Define $\fnogiant(n,m) \eqdef \E(\text{\MGF}(G(n,m))$.
\begin{theorem}
\label{bfwplus}
With $t=t(c)< 1$ defined by $te^{-t}=2ce^{-2c}$,
$\fnogiant(n,cn) = rcn$
when
$\frac{t^2}{4c}+1= \frac{t}{2c}+ c r$.
\end{theorem}
The theorem is proved as in~\cite{BFW02}
(modifying their Lemma~1 to allow values $c>2$ by replacing a
$(\log n)/6$ with $(\log n)/(6\log c)$).

Is there another \maxsc subgraph problem, then, which does behave
like \maxtsat?
Going back to the branching process for a random graph ---
the source of the intuitive connection between the graph and \sat problems ---
it is also easy to check that \whp a graph has few cycles
when the branching process is subcritical,
and many cycles when it is supercritical.
So perhaps we should consider the size of maximum cycle-free subgraph.
But this is by definition a forest, which may have at most $n-1$ edges,
again a $1/c$ fraction, not a fixed constant fraction as for \maxtsat.

In a 2-\sat formula, obstructions to satisfiability come not from
cycles of implications $X \implies \cdots \implies X$,
but only from those with $X \implies \cdots \implies \Xbar$.
By a very vague analogy, then, perhaps on the graph side we should
seek not a subgraph which is entirely cycle-free,
but just one which is free of \emph{odd} cycles:
a bipartite subgraph.
The size of a largest bipartite subgraph $H$ of $G$ is
by definition, and more familiarly,
the size of a maximum cut of $G$.
Here, finally, we share with \maxtsat that we may keep a constant
fraction of the input structure:
for a random graph (indeed any graph) $G$ of size $m$,
$\maxcut(G) \geq m/2$, since a random cut achieves this expectation.

\subsection{MAX CUT}
In addition to the fact that just as a maximum assignment satisfies
at least $3/4$ths the clauses of any formula,
a maximum cut cuts at least $1/2$ the edges of a graph,
there are other commonalities.

\maxcut, like \maxtsat,
is a constraint satisfaction problem (\csp).
With each vertex $v$ we associate a boolean variable representing
the partition to which $v$ belongs,
and with each edge $\set{u,v}$
we associate a ``cut constraint'' $(u \oplus v)$,
these \xor constraints replacing 2-\sat's disjunctions.
Like decision 2-\sat,
the problem of whether a graph is perfectly cuttable (bipartite)
is solvable in essentially linear time.
In further analogy with \maxtsat, \maxcut is \np-hard,
trivially $\frac12$-approximable,
0.878-approximable~\cite{GW2} by semidefinite programming,
and not better than $16/17$-approximable \cite{SSTW00}
in polynomial time, unless P=NP.

The methods we have applied to random \maxtsat are equally
applicable to \maxcut, and yield analogous results.
Because it is easier to work with random graphs than random formulas,
and more is known about them, our results for \maxcut are in some
respects stronger than those for \maxtsat.

When we work in the \gnp model we will take $p=2 c/n$,
and in the \gnm model, $m=\lfloor c n \rfloor$,
so that in both cases the phase transition occurs at~$c=1/2$.
We now state our main results.

\subsection{Results}

\begin{theorem}
\label{CutSubcritical}
For $c = 1/2 - \e(n)$, with $n^{-1/3} \ll \e(n) < 1/2$,
$\fc(n,\cn) = \cn - \Theta(\ln(1/\e)) + \Theta(1)$.
\end{theorem}
In particular, for small constants $\e$ this gap of $\Theta(\ln(1/\e))$
--- which for a fixed $\e$ is $\Theta(1)$ ---
contrasts with the gap of
$\Theta(1/n)$ for \maxtsat (Theorem~\ref{subthreshold}).
But here too there is a phase transition,
in that for $c>1/2$ the gap jumps to $\Theta(n)$,
per Theorem~\ref{CutLowDens}.

\can{\cuthigh}{%
\begin{theorem}\label{CutHighDens}
For $c$ large,
$\left( \frac12 c + \sqrt{c} \cdot \sqrt{8/(9 \pi)} \right) n
\leqa \fc(n,cn) \leqa
\left( \frac12 c + \sqrt{c} \sqrt{\ln(2)/2} \right) n$.
\end{theorem}
}
\uncan{\cuthigh}

\noindent
The values of $ \sqrt{8/(9 \pi)}$ and $\sqrt{\ln(2)/2}$ are
approximately $0.531922$ and $0.588704$, respectively.
The upper bound
was previously obtained in~\cite{Bertoni97}.

\can{\cutlow}{%
\begin{theorem} \label{CutLowDens}
For any fixed $\e>0$,
$(\frac12 + \e - \frac{16}{3}(\e^3)) n
\leqa \fc(n,(1/2+\e)n) \leqa
(\frac12 + \e - \Omega(\e^3/ \ln(1/ \e))) n$.
\end{theorem}
}
\uncan{\cutlow}
\noindent
The upper bound's $\e^3/\ln(1/\e)$
can probably be replaced by~$\e^3$,
just as we suspect it can be for Theorem~\ref{MainLowDens}.
This presumption is largely based on the next ``scaling window'' result.

\begin{theorem}
\label{cutscaling}
For any function $\e = \e(n)$ with $n^{-1/3} \ll \e(n) \ll 1$,
$\fc(n,(1/2+\e)n) =
(\frac12 + \e - \Theta(\e^3)) n$.
\end{theorem}
That the theorem misses out the extremes $\e=\Theta(n^{-1/3})$
and $\e=\Theta(1)$ that are perhaps of greater interest than
the mid-range is a direct carryover from the standard results
on random graphs on which we based our proof is based;
it is likely that other established results for random graphs
could complete the picture.

Before proceeding, we remark
that bipartiteness is of course the same as 2-colorability, and
it is sometimes convenient to speak of coloring vertices
black or white,
rather than placing them in the left or right part of a partition,
with properly colored edges (with one black and one white endpoint)
corresponding to cut edges;
these two ways of speaking are of course mathematically identical.

\subsection{Subcritical MAX CUT}

\thmpf{CutSubcritical}{Proof.}
For notational convenience we work in the \gnp model,
$G = G(n, (1-2\e)/n)$,
but the proof follows identically for the \gnm model.

Tree components of $G$ can be cut perfectly;
each unicyclic component can be cut for all but 1 edge;
and complex components,
where more edges must go uncut
but which with high probability are absent from~$G$,
contribute negligibly.
That is, $\E(\text{\#uncut edges}) = (1-o(1)) \E(\text{\#cycles in $G$})$.
Since the number of potential $k$-cycles is $(n)_k/(2k)$,
where $(n)_k = n(n-1)\cdots(n-k+1)$ denotes falling factorial,
using $(n)_k = n^k \exp(-k^2/(2n) - O(k/n + k^3/n^2))$
(see~\cite[eq~(5.5)]{JLR}),
\begin{align*}
\E(\text{\#cycles in $G$})
 &= \sum_{k=3}^n \frac{(n)_k}{2k} (c/n)^k
 \\ &= \sum \frac{1}{2k} c^k \exp(-k^2/(2n) - O(k/n + k^3/n^2)) .
\end{align*}
Because of the $c^k$, up to constant factors we need consider the
sum only up to $k \leq 1/\e$ (recalling $c=1-2\e$),
and since $\e \gg n^{-1/3}$,
this makes the entire final exponential term negligibly close to~1.
Thus
\begin{align*}
\E(\text{\#cycles in $G$})
 &= \Theta(1) \sum_{k = 3}^{\infty} c^k/(2k)
 \\ &= \Theta(1) (-\frac12 \ln(1-c)) -  \Theta(1)
 \\ &= \Theta(1) \ln(1/(2\e)) - \Theta(1) ,
\end{align*}
where the final $\Theta(1)$ term lies between 0 and $3/2$.
\qed

\subsection{High-density random MAX CUT}

\thmpf{CutHighDens}{Proof.}
For the upper bound, we apply a first-moment argument identical to that used in the
proof of Theorem~\ref{MainHighDens}.
The probability that there exists a (maximal) bipartite spanning subgraph
of size $\geq (1-r)cn$ is
$P \leqa 2^n \binom{cn}{rcn} (1/2)^{(1-r)cn} (1/2)^{rcn}$,
for $\frac{1}{cn} \ln P \leqa \ln 2/c -\logis{r}-\logis{(1-r)}-\ln 2$.
Substituting $r=1/2-\e$ gives
$\frac{1}{cn} \ln P \leqa \ln 2 / c - 2 \e^2$,
so if $\e > \sqrt{\ln(2)/(2c)}$ then $P \ra 0$.

For the lower bound, color the vertices in random sequence.
When $xn$ vertices have been colored, with $x=\Theta(1)$,
since $c$ is large,
the next vertex is a.s.\ adjacent to a.e.\ $2 c x$
of the colored ones.
In the worst case, the colored vertices are half black and half white;
coloring the new vertex oppositely to the majority beats $cx$
(in expectation) by
$\E(|B(2cx, 1/2)  - cx|) \asymp \E(|N(0, cx)|) = \sqrt{2cx/\pi}$.
Integrating over $x$ from 0 to 1 gives $n \sqrt{2c/\pi} \cdot \frac23$
more properly colored edges than the naive $\frac12 c n$.
\qed

\subsection{Low-density random MAX CUT}

The following fact follows from small-$\e$ asymptotics of classical
random graph results;
see, e.g.,~\cite[VII.5, Theorem~17]{ModernBollobas98}.

\begin{claim}
For $\e>0$, a random graph $G(n,(1/2+\e)n)$ a.s.\ has a giant
component of size $(4 \e +o(\e)) n$.
\end{claim}

\begin{proof}
It is well known
(see, e.g.,~\cite[VII.5, Theorem~17]{ModernBollobas98}) that
for an arbitrarily slowly growing function
$\w(n)$, a.s., the size $L^{(1)}(G)$ of the giant component
satisfies $|L^{(1)}(G) - \gamma n| \leq \w(n) n^{1/2}$ where
$0<\gamma<1$ is the unique solution of $e^{-2c \gamma} = 1-\gamma$.
(We have $2c$ where \cite{ModernBollobas98} has $c$ because
we use $cn$ edges where it uses average degree $c$.)
Take the asymptotic approximation when $c=1/2+\e$.
\end{proof}

\begin{claim}
\label{bipartite}
The probability that a random graph $G(n,(1/2+\e)n)$ is bipartite,
conditioned on the existence of a component of
size $\Theta(\e n)$ created by the ``first'' $(1/2+\e/2)n$ edges,
is $\exp(-\Omega(\e^3 n))$.
\end{claim}

\begin{proof}
If the presumed giant component is not bipartite, we are done.
If it is, by connectivity, it has a unique bipartition;
let the sizes of the parts be $n_1$ and $n_2$.
Each of the remaining $\e n/2$ edges
has both endpoints in the giant component w.p.~$\Theta(\e^2)$,
so there are $\Theta(\e^3 n)$ of these,
w.p.\ $1-\exp(-\Omega(\e^3 n))$.
The probability that each such edge preserves bipartiteness
is $(2n_1n_2)/(n_1+n_2)^2\leq 1/2$;
over the $\Theta(\e^3 n)$ independent edges it is $\exp(-\Omega(\e^3 n))$.
\end{proof}

\thmpf{CutLowDens}{Proof.}
For the upper bound, the first-moment method is applied exactly as in the proof
of Theorem~\ref{MainLowDens}.
We use the preceding Claim, and replace its $\Omega$ with an $\a_0$
for definiteness.
With $c=(1/2+\e)$, then,
the probability that deleting any $k \leq rcn$ edges
can leave a bipartite subgraph is
$P \leq \sum_{k=0}^{rcn} \binom{cn}{k} \exp(-\a_0 (\e-k/n)^3)$.
This is just as in inequality~\eqref{logbound},
so here again we conclude that $r \geqa \a_0 \e^3 / \ln(1/\e)$.

The proof of the lower bound is algorithmic, and in direct analogy to
that of
Theorem~\ref{MainLowDens}.
Think of a graph edge neither of whose vertices
has yet been colored as a ``2-clause'',
an edge one of whose vertices has been colored as a ``unit clause''
implying the opposite color for the remaining vertex,
an edge whose two vertices
have been colored alike as an ``unsatisfied clause'',
and an edge whose
two vertices have been colored oppositely as a ``satisfied clause''.
Terminate if there are no unit clauses nor 2-clauses.
If there are no unit clauses,
randomly color a random vertex from a random edge.
If there are unit clauses,
choose one at random and color its vertex satisfyingly.

Note that when a $\d$ fraction of the vertices have been colored,
$(1-\d)n$ vertices remain uncolored (unfixed variables), and a.s.\
a.e.\ $(1/2+\eps)n \cdot (1-\d)^2$ 2-clauses remain.
Each time a unit-clause variable is set,
each 2-clause has probability
$2/((1-\d)n)$ of generating a unit clause

Thus the expected number of 2-clauses becoming unit clauses is $2
[(1/2+\e)n (1-\d)^2] / [(1-\d) n] \approx 1+ 2 \e - \d$, while the
number of unit clauses eliminated (satisfied or unsatisfied) is at
least~1. Thus the expected increase per step in the number of unit
clauses is at most $2\e-\d$. As in the proof of Theorem~\ref{MainLowDens},
over the first $4\e n$ steps, the expected number of unit clauses
is bounded by an inverted parabola of base $4 \e n$ and height $2 \e^2 n$.
Improperly colored edges result only from violated unit clauses,
and the expected number of these in the first $4\e n$ steps is
$\leq \frac23 \cdot 4\e n \cdot 2\e^2n / n$ $\leq \frac{16}{3}
\e^3$. By step $4\e n$ there are no unit clauses, and the number
of 2-clauses divided by the number of unset variables is a.s.\ a.e.\
$[(1/2+\e)n \cdot (1-4\e)^2] / [n \cdot (1-4\e)] = 1/2-\e$.
This is a sparse random graph,
which by Theorem~\ref{CutSubcritical} can be colored to violate
just $\Theta(1)$ edges.

\textit{In toto}, all but $\leqa (\frac{16}{3} \e^3 n + \Theta(1))$
edges are properly colored.
\qed

\subsection{Scaling window}
\label{sec:maxscaling}

The proof of Theorem~\ref{cutscaling} follows rather easily from
standard --- but relatively recent, and lovely --- facts about
the kernel of a random graph.
The following summary of the relevant facts,
which we present informally, is distilled from~\cite[Sec.~5.4]{JLR}.

\renewcommand{\w}{\omega}
First, if $i \gg n^{2/3}$, then
the number of vertices of $G(n, n/2+i)$ belonging
to unicyclic components is asymptotically almost surely
$\Theta(n^2/i^2)$.
Consider the components of a graph $G$ which are trees,
unicyclic, or complex.
In the supercritical phase with $n^{-1/3} \ll \e \ll 1$,
a random graph $G(n,(1/2+\e)n)$ consists of tree components,
unicyclic components, and no complex component other than
a single ``giant component''.
The expected number of vertices in the cycles of the unicyclic components
is of order $1/\e$.
The giant component's 2-core has order $(1+o(1)) 8 \e^2 n$,
and is obtained as a random subdivision of the edges of
a ``kernel'', which is a random cubic graph on
$(1+o(1)) \frac{32}{3} \e^3 n$ vertices.

\thmpf{cutscaling}{Proof.}
We consider which edges of $G$ it may be impossible to cut.
Every edge in the tree components of $G = G(n,(1/2+\e)n)$
can of course be cut.
For each unicyclic component, at most 1 edge must go uncut
(if the cycle is odd).
By the symmetry rule (see for example~\cite[Theorem~5.24]{JLR}),
the number of unicyclic components for $G(n,(1/2+\e n))$
is essentially the same as for $G(n,(1/2-\e n))$,
which by Theorem~\ref{CutSubcritical} is
only $O(\ln(1/\e))$.

The dominant contribution will come from the giant component.
Edges which are not in its 2-core can of course all be cut,
even after a partition of the 2-core has been decided.
Moreover, an optimal partition of the 2-core is essentially
decided by a partition of the vertices of the ``kernel'',
which is the 2-core where each path whose internal vertices
are all of degree~2 is replaced by a single edge.
(See \cite[Chap.~5.4]{JLR} for more on the giant component,
its core, and its kernel.)
For any cut of the kernel,
each 2-core path corresponding to a kernel edge
can be partitioned either perfectly or with one edge uncut,
depending on the parity of the path's length
and whether its endpoints are on the same side or opposite
sides of the kernel's cut.
Equivalently, a
kernel edge whose 2-core path is of odd length
imposes a ``cut'' constraint on its endpoints,
while a kernel edge whose 2-core path is of even length
imposes an ``uncut'' constraint on its endpoints;
the number of these constraints violated by a cut of the
kernel vertices is equal to
the number of original cut constraints violated by an
optimal extension of the same cut to all the 2-core vertices
(and indeed to all the giant-component vertices).

Since each kernel edge is randomly subdivided,
on average into $3/(4\e)$ 2-core edges,
the parities of the kernel edges are almost perfectly random
(with the probability of either parity approaching $1/2$
as $\e$ approaches~$0$).
For our purposes it suffices that either parity occurs
with probability at most some absolute constant $p_0 < 1$,
and using this we show that at least some constant fraction $\b_0$
of the approximately $16 \e^3 n$ edge constraints must be violated.

Fix a spanning tree $T$ of the kernel~$K$, whose order
we will write as $N$ (expecting $N \approx \frac{32}{3} \e^3 n$).
Let $K$ subsume not only the graph but also the edge parities,
so that it is an instance of the generalized (cut/uncut) \maxcut problem.
If it is possible to violate precisely a fraction $\b < \b_0$
of $K$'s constraints
then reversing precisely those constraints gives a perfectly
satisfiable cut/uncut constraint problem instance~$K'$.

Fixing the ``side'' of any one vertex,
the $N-1$ constraints from the spanning tree $T$
imply the rest of the cut,
which must then satisfy the remaining $\frac{1}{2}|N|+1$ constraints.
Viewing the parities of the spanning tree edges as arbitrary,
and the remaining edges as independent random variables,
the probability that the randomly chosen kernel edges
satisfy each of these constraints is at most $p_0^{\frac{1}{2}N+1}$.
The number of choices of $t < \beta_0 \frac{3}{2}N$ edges
to dissatisfy is ${{\frac{3}{2}N} \choose t}$.
We guarantee an exponentially small probability of success
by selecting $\beta_0$ to satisfy:
\begin{align*}
\sum_{t < \beta_0 \frac{3}{2}N }
  {{\frac{3}{2}N} \choose t}  p_0^{\frac{1}{2}N+1}
  & \ll 1
\\
\frac{3}{2}N  H(\beta_0) +   \frac{1}{2}N  \ln ( p_0)
 & < 0
\\
H(\beta_0)
 &<   \frac{1}{3}  \ln ( 1/p_0) ,
\end{align*}
where $H$ is the entropy function $H(x) = x \ln(x) - (1-x) \ln(1-x)$.
In particular, in the case of interest where $\e \ra 0$,
$p_0 \ra 1/2$ and $\b_0 \ra H^{-1}(1/3) \approx 0.896$.
Recapitulating, we must dissatisfy $\b_0 N$ kernel constraints,
$=(32 \b_0 / 3) \e^3 n$ constraints of $G$.
The expected $O(\ln(1/\e))$ uncut edges from unicyclic components
are negligible by comparison, so in all
$\Theta(\e^3 n)$ edges of $G$ go uncut.

\qed

\section{Conclusions and open problems}
We have presented a road map for \maxtsat and \maxcut in a random setting,
establishing that there is a phase transition,
and deriving asymptotics below the critical value,
for constants slightly above the critical value and
in the scaling window around it,
and for larger constants.

For constant densities slightly above threshold there is a
logarithmic gap between our lower and upper bounds;
we need to confirm that the $\ln(1/\e)$ factors are extraneous.
In the other cases, our bounds are only separated by a constant.
However, in light of the exact result of \cite{BFW02}
for the size of a maximum subgraph which has no giant component,
it would be wonderful to get the \emph{exact} asymptotics
of $f(n,cn)/(cn)$.

Whether $f(n,cn)/(cn)$ tends to a limit in $n$
(see Conjecture~\ref{conjecture:limit}) is to our minds
a prime open problem in this area,
and is not only in some sense analogous to the
satisfiability threshold conjecture,
but may also be directly connected with it
(see Conjecture~\ref{conj2}),
another important question.

A question similar in spirit to Conjecture~\ref{conjecture:limit}
was considered in~\cite{Gamarnik},
which defines a certain
linear-programming relaxation of \maxtsat.
An instance in characterized by its ``distance to feasibility''~$O$,
with $O(n,cn)$ the corresponding
random variable
for a random instance.
It is shown that for every $c>0$,
$O(n,cn)/(cn)$ almost surely converges to a limit.
The result is established using powerful local weak convergence methods
\cite{AldousAssignmentI,AldousAssignmentII,AldousSteele}.
It remains to be seen whether these methods are applicable to
random maximum constraint satisfaction problems,
including \maxtsat and \maxcut.

\section*{Acknowledgments}
Dr~Sorkin is very grateful to Svante Janson for pointing out
the beautiful properties of the kernel of a random graph
and its relevance, and for other helpful technical discussions.
Dr~Sorkin is also grateful to Mark Jerrum, Martin Dyer,
and Peter Winkler for organizing the excellent Newton Institute
workshop which fostered these exchanges.


\newcommand{\etalchar}[1]{$^{#1}$}
\providecommand{\bysame}{\leavevmode\hbox
to3em{\hrulefill}\thinspace}
\providecommand{\MR}{\relax\ifhmode\unskip\space\fi MR }
\providecommand{\MRhref}[2]{%
  \href{http://www.ams.org/mathscinet-getitem?mr=#1}{#2}
} \providecommand{\href}[2]{#2}

\end{document}